%% file: linear-patterns.tex
\documentclass[11pt]{article}
\usepackage[margin=1in]{geometry}
\usepackage{thm-restate}
\usepackage{amsmath,amssymb,amsthm}
\usepackage[hidelinks]{hyperref}
\usepackage{graphicx}
\usepackage{tikz}
\usepackage{arydshln}
\usepackage{stmaryrd}

\usepackage[capitalise]{cleveref}

\hypersetup{pdftitle={Linear Extremal Functions for Forbidden Matrices},
pdfauthor={Jesse Geneson; Seth Pettie}}

\newtheorem{theorem}{Theorem}[section]
\newtheorem{lemma}[theorem]{Lemma}
\newtheorem{corollary}[theorem]{Corollary}

\newtheorem{conjecture}[theorem]{Conjecture}

\Crefname{conjecture}{Conjecture}{Conjectures}

\theoremstyle{remark}
\newtheorem{remark}[theorem]{Remark}
\newtheorem{definition}[theorem]{Definition}

\newcommand{\ex}{\operatorname{ex}}

\newcommand{\bu}{\ensuremath{\bullet}}
\long\def\restateable#1#2#3{\begin{restatable}{#1}{#2}#3\end{restatable}}

\newcommand{\nf}{\acute{n}}
\newcommand{\nl}{\grave{n}}
\newcommand{\nm}{\bar{n}}
\newcommand{\Uf}{\acute{U}}
\newcommand{\Ul}{\grave{U}}
\newcommand{\Um}{\bar{U}}
\newcommand{\Sf}{\acute{S}}
\newcommand{\Sl}{\grave{S}}
\newcommand{\Sm}{\bar{S}}

\newcommand{\calS}{\mathcal{S}}
\newcommand{\calSf}{\acute{\mathcal{S}}}
\newcommand{\calSl}{\grave{\mathcal{S}}}
\newcommand{\calSm}{\bar{\mathcal{S}}}

\newcommand{\dbl}{\operatorname{dbl}}
\newcommand{\FH}{\operatorname{FH}}
\newcommand{\SW}{\operatorname{SW}}
\newcommand{\Ge}{\operatorname{Ge}}

\newcommand{\first}{\text{first}}
\newcommand{\last}{\text{last}}
\newcommand{\lo}{\text{lo}}
\newcommand{\hi}{\text{hi}}
\newcommand{\Span}{\text{span}}

\newcommand{\Top}{\text{top}}
\newcommand{\Bottom}{\text{bot}}

\newcommand{\hcm}[1][1]{\hspace*{#1 cm}}
\newcommand{\rb}[2]{\raisebox{#1 mm}[0mm][0mm]{#2}}
\newcommand{\istrut}[2][0]{\rule[- #1 mm]{0mm}{#1 mm}\rule{0mm}{#2 mm}}
\newcommand{\zero}[1]{\raisebox{0pt}[0pt][0pt]{\makebox[0pt]{$#1$}}}
\newcommand{\ceil}[1]{\left\lceil #1 \right\rceil}

\newcommand{\Furedi}{F\"{u}redi}

\newcommand{\Turan}{Tur\'{a}n}

\newcommand{\ignore}[1]{}

\input{patbits}

\title{Classification of Forbidden 0--1 Patterns\\ by Their Extremal Function\thanks{Supported by NSF Grant CCF-2446604.}}

\author{Jesse Geneson\\San Jose State University \and Seth Pettie\\University of Michigan}
\date{}
\begin{document}
\maketitle

\begin{abstract}
This paper concerns the maximum weight of $n\times n$ 0--1 matrices 
that avoid some pattern $P$, denoted $\ex(n,P)$.
\Furedi{} and Hajnal~\cite{FurediH92} initiated a research program of
classifying all patterns according to their extremal function, but 
to date, only the weight-4 patterns are completely understood.

\medskip
In this paper we make progress on several fronts of the classification problem.  A pattern is \emph{light} if its columns have weight 1 and \emph{linear} if $\ex(n,P)=O(n)$.

\begin{itemize}
    \item We characterize the set of weight-5 linear light patterns, and the set of weight-5 linear patterns more generally.  Along the way we identify several new infinite classes of linear patterns, which support the Pettie--Tardos conjecture on characterizing all linear light patterns.
    \item It was known that there \emph{exist} infinitely many minimally nonlinear patterns.  We identify a \emph{specific} infinite set of minimally nonlinear patterns. 
    We also prove that there are $2^{\Omega(r)}$ 
    weight-$r$ minimally nonlinear patterns, and that every minimally nonlinear pattern has extremal function $O(n\log n)$, with the exception of the 4-cycle.
    \item We classify many of the outstanding weight-5 patterns according to their extremal function, leaving only a 
    small set remaining.
    The newly classified weight-5 patterns have extremal function $\Theta(n)$ or $\Theta(n\alpha(n))$.  We also exhibit a large number of weight-6 patterns with extremal function 
    $\Theta(n\log n/\log\log n)$, which is notable
    in that this function 
    had not been encountered before for
    any \emph{single} pattern $P$, though it was known to be the correct extremal function for a \emph{pair} of forbidden patterns.
\end{itemize}

\end{abstract}

\section{Introduction}\label{sec:introduction}

The theory of pattern avoiding 0--1 matrices has many applications in geometry \cite{FurediH92,PachTardos06}, enumerative combinatorics~\cite{MarcusT04,Klazar00,Cibulka09}, 
structural graph theory~\cite{GuillemotM14,BonnetNMST24,BonnetGMSTT24},
and the analysis of data structures~\cite{Pettie10a}, 
particularly special cases of the dynamic optimality conjecture~\cite{ST85,DemaineHIKP09,ChalermsookGKMS15,ChalermsookGKMS15b,ChalermsookGJAPY23,ChalermsookPY24}.  
We say that a pattern $P\in \{0,1\}^{k\times l}$ \emph{appears in} or \emph{occurs in} $A\in\{0,1\}^{n\times m}$, written $P\prec A$, if there is a $k\times l$ submatrix of $A$ that dominates $P$.  In other words, a 1 in $P$ must match a 1 in $A$, but a 0 in $P$ matches anything. We say $A$ is \emph{$P$-free} if $P\nprec A$.
Define $|A|$ to be the \emph{weight} or number of 1s in $A$.
The main object of study is the extremal function
\[
 \ex(n,m,\mathcal P)=
 \max\{|A| \;:\; A\in\{0,1\}^{n\times m} \text{ and } \ P\nprec A
                         \text{ for every }P\in\mathcal P\}.
\]
This notation permits rectangular $A$ and a set of forbidden patterns,
but most of the time we shall focus on square matrices avoiding a single forbidden pattern.  The simplified extremal functions are:
\begin{align*}
    \ex(n,m,P) &= \ex(n,m,\{P\}),\\
    \ex(n,P) &= \ex(n,n,P).
\end{align*}

This subject is usually regarded as an \emph{ordered} bipartite version of \Turan{} theory~\cite{FurediS13}, where vertices on each side of the bipartition are ordered.
Whereas \Turan{} extremal functions are always $O(n)$ or $n^{1+\Omega(1)}$,
the extremal functions for ordered 0--1 
patterns are more varied~\cite{FurediH92,Pettie-SoCG11,Pettie-FH11,KorandiTTW19,PettieT24,PettieT25}.  
All known \emph{sharp} bounds that are $n^{1+o(1)}$ are asymptotic to one of the following.
\[
n,
\; n\alpha(n),
\; n2^{\alpha(n)},
\; n\alpha(n)2^{\alpha(n)},
\; n2^{(1+o(1))\alpha(n)^t/t!},
\; n\log n,
\; n(\log n/\log\log n)^t,
\; n2^{\Theta(\sqrt{\log n})},
\]
where $t\geq 2$ is an integer and $\alpha(n)$ is the slowly growing 
inverse-Ackermann function.  Some functions arise that are not known to be sharp, e.g., $n\alpha^2(n)$~\cite{Pettie-GenDS11}, 
$O(n\log^t n)$ and 
$n2^{(\log n)^{1-1/(t+1)}}$,
for $t\geq 2$~\cite{PachTardos06,KorandiTTW19}, 
and $n^{1+O(1/\log\log n)}$.  
The last, due to 
Gishboliner and Li~\cite{GishbolinerL26}, applies to \emph{all} acyclic patterns, positively resolving the weak Pettie--Tardos conjecture~\cite[Conjecture 4.1]{PettieT25}.

\medskip

A pattern is \emph{light} if every column contains exactly one 1.
Matrices avoiding light patterns are closely related to generalized 
Davenport--Schinzel sequences~\cite{Klazar92,Keszegh09,Pettie-SoCG11,Pettie-GenDS11,Pettie15-SIDMA}.
It is known~\cite{Keszegh09} 
that every light $P$
has $\ex(n,P) = O(nf(\alpha(n)))$, 
where $f$ is one of the functions of 
$\alpha(n)$ listed above.

\subsection{The Classification Problem}

Over 35 years ago, \Furedi{} and Hajnal~\cite{FurediH92} initiated a research program to \emph{classify} all 0--1 patterns according to their extremal function.  They classified all weight-4 patterns except for two, which Tardos later resolved~\cite{Tardos05}.  We now know that
all \emph{light} weight-4 patterns are linear except for 
$Q_3,Q_3'$, which are $\Theta(n\alpha(n))$.
\[
 Q_3=\patword{1212},\qquad
 Q'_3=\patword{1323}.
\]
Order has no effect on the all-1 
$2\times 2$ pattern $C_4$,
which has $\ex(n,C_4)=\Theta(n^{1.5})$~ \cite{KovariST54,Reiman58}.
The only other superlinear weight-4 patterns $R_1,R_2,R_3$
have extremal function 
$\Theta(n\log n)$~\cite{Tardos05,FurediH92}.
\begin{align*}
C_4 &= \patbits{11/11},
&R_1 &= \patbits{110/101},
&R_2 &= \patbits{100/010/101},
&R_3 &= \patbits{110/100/001}.
\end{align*}

\Furedi{} and Hajnal asked for a characterization of all linear patterns (those with $\ex(n,P) = O(n)$) and specifically conjectured that $\ex(n,P)=O(n)$ whenever $P$ is a permutation matrix.  This conjecture was proved by Marcus and Tardos~\cite{MarcusT04},
which also implied the Stanley-Wilf conjecture via a prior reduction of Klazar~\cite{Klazar00}.
Previously known linear families include Geneson's~\cite{Geneson09} 
``double'' permutation patterns (generalizing~\cite{MarcusT04}) and Keszegh's~\cite{Keszegh09,Pettie-SoCG11} 
monotone patterns.  An example of each is given below.\footnote{It is easy to see that in a light pattern, repeating the first or last column more than once, or any interior column more than twice, does not affect the asymptotics of the extremal function.  However, there is no theorem to the effect that ``doubling'' a light pattern has no effect on its extremal function.  See~\cite{Pettie-GenDS11,Pettie15-SIDMA}.}

\begin{align*}
    \patword{52233114}, & & & \patword{5443322112233445}.
\end{align*}

Examples of linear patterns discovered by
Tardos~\cite{Tardos05}, Fulek~\cite{Fulek09},
Pettie~\cite{Pettie-GenDS11}, and Geneson~\cite{Geneson26}
are shown below.
\begin{align*}
    L_{T} &= \patbits{0100/1000/0001/1010},
    &&& L_{F} &= \patbits{00100/10001/01110},\\
    L_{P} &= \patbits{0000010/0110000/0001101/1000000},
    &&& L_{G} &= \patbits{10011/000001/011000}.
\end{align*}

Keszegh~\cite{Keszegh09} and Geneson~\cite{Geneson09} ruled out the simplest characterization of the class of linear patterns, namely a \emph{finite} 
list of minimally
nonlinear patterns with respect to containment ($\prec$).  
They proved that the set of minimally nonlinear patterns is infinite, but could not identify a \emph{particular} 
infinite set of minimally nonlinear patterns.
Nonetheless, all of Keszegh's patterns shown to exist 
have extremal function $\Theta(n\log n)$.  In particular they are not light.  Pettie and Tardos~\cite{PettieT25} recently conjectured that the class of linear light patterns \emph{is}, in fact, 
characterized by a finite list of minimally nonlinear patterns.

\begin{conjecture}[{Pettie and Tardos~\cite[Conjecture 4.2]{PettieT25}}]\label{conj:PT}
Define $\mathcal{Q}$ to be the set of patterns containing $Q_3,Q_3'$ and 
their reflections.  If $P$ is a light pattern, then $\ex(n,P)=O(n)$ if and only if $P$ is $Q$-free for all $Q\in\mathcal{Q}$.
\end{conjecture}

\subsection{Pattern Transformations}

Upper and lower bounds on $\ex(n,P)$ can often be deduced via \emph{pattern transformations}.

\paragraph{Trivial transformations.} \Furedi{} and Hajnal observed that the extremal function $\ex(n,P)$ is invariant under rotation and reflection of $P$,
and that $P\prec P'$ implies $\ex(n,P) \leq \ex(n,P')$.  Moreover, suppose
$P\in\{0,1\}^{k\times l}$ is such that for some row $r$, $P(r,l)=1$.  
If $P'\in\{0,1\}^{k\times (l+1)}$
is obtained from $P$ by appending a weight-1 column with $P'(r,l+1)=1$,
then $\ex(n,P') \leq n + \ex(n,P)$.
For example, these transformations alone imply the following sequence has 
non-increasing extremal functions, asymptotically.
\[
 \patbits{11/11}\ \longrightarrow\
 \patbits{10/01/11}\ \longrightarrow\
 \patbits{10/01/10/01}\ \longrightarrow\
 \patbits{100/001/010/001}\ \longrightarrow\
 \patbits{1000/0010/0100/0001}.
\]

\paragraph{Keszegh's Gluing Operation.}

Keszegh~\cite{Keszegh09} proved that whenever $A$ has a 1 in the SE corner and $B$ has a 1 in the NW corner, we can identify these 1s, so that the resulting pattern $A\oplus B$ satisfies
$\ex(n,A\oplus B) \leq \ex(n,A) + \ex(n,B)$.
\begin{align*}
A\oplus B &= 
\left(\,\begin{array}{ccccccc}
\cline{1-4}
\multicolumn{1}{|c}{ } &&& \multicolumn{1}{c|}{} \\
\multicolumn{1}{|c}{ } &\zero{\hcm[.25]{A}} && \multicolumn{1}{c|}{} &&& \\\cline{4-7}
\multicolumn{1}{|c}{ } &&\zero{\hcm[.7]{\bu}}& \multicolumn{1}{|c|}{} &&& \multicolumn{1}{c|}{}\\\cline{1-4}
&&& \multicolumn{1}{|c}{} &&\zero{\hcm[-.2]{B}}& \multicolumn{1}{c|}{}\\
&&& \multicolumn{1}{|c}{} &&& \multicolumn{1}{c|}{}\\\cline{4-7}
\end{array}\,\right).
\end{align*}

\paragraph{Pach--Tardos Column Addition.}

Pach and Tardos~\cite{PachTardos06}
showed that in two circumstances, 
adding a weight-1 column can increase the extremal function by at most a 
$\log n$ factor.
\begin{lemma}[\cite{PachTardos06}]\label{lem:PachTardos-reductions}
Suppose $P'$ is obtained from $P$ 
by appending a weight-1 column.
Then $\ex(n,P') = O(n+\ex(n,P)\log n)$.
Suppose $P$ has two consecutive columns 
$j,j+1$ and rows $i,i'$ with 
$P(i,j)=P(i',j)=P(i,j+1)=1$. 
If we obtain $P'$ by splicing in a 
weight-1 column ``$j+1/2$'' with 
$P'(i',j+1/2)=1$, then
$\ex(n,P') = O(\ex(n,P)\log n)$.

\begin{align*}
    P' &= \left(\,\begin{array}{|ccc|c}
    \cline{1-3}&&&\\
    &P&&\bu\\
    &&&\\\cline{1-3}\end{array}\,\right),
    &&&
    P' &= \left(\,\begin{array}{|ccc|c|ccc|}
    \cline{1-3}\cline{5-7}&&\bu&&\bu&&\\
    &P&&&&&\\
    &&\bu&\bu&&&\\\cline{1-3}\cline{5-7}\end{array}\,\right).
\end{align*}
\end{lemma}

\paragraph{Pettie's Grafting Operation.}
Suppose $P_0$ is a pattern with 1s in the SW and SE corners, and $P_1$ is a pattern with two 1s in the top row in consecutive columns $j,j+1$.
Obtain $P_0\otimes P_1$ by identifying the 
SW 1 of $P_0$ with the 1 of $P_1$ in column $j$, and inserting the remaining columns of $P_0$ strictly before column $j+1$ of $P_1$.
This is referred to as \emph{grafting} $P_0$ onto $P_1$.
See the diagram below.

\[
\begin{array}{ccccccc}
&\hcm&\zero{P_0}&\hcm[.4]&&&\hcm\\\cline{3-5}
&&\multicolumn{1}{|c}{} && \multicolumn{1}{c|}{}\\
&&\multicolumn{1}{|c}{} && \multicolumn{1}{c|}{} &&\\
\zero{P_1}&&\multicolumn{1}{|c}{} &&\multicolumn{1}{c|}{}\\\cline{1-3}\cline{6-7}
\multicolumn{1}{|c}{}&&\multicolumn{1}{|c|}{\bu} &&\multicolumn{1}{c|}{\bu} & \bu&\multicolumn{1}{c|}{}\\\cline{3-5}
\multicolumn{1}{|c}{}&&\multicolumn{1}{c|}{} &&\multicolumn{1}{c|}{} & &\multicolumn{1}{c|}{}\\
\multicolumn{1}{|c}{}&  &\multicolumn{1}{c|}{} &&\multicolumn{1}{c|}{} & &\multicolumn{1}{c|}{}\\
\multicolumn{1}{|c}{}&&\multicolumn{1}{c|}{} &&\multicolumn{1}{c|}{} & &\multicolumn{1}{c|}{}\\\cline{1-3}\cline{6-7}
\end{array}
\]

\begin{definition}[Legality]\label{def:legal}
A pattern $P$ is \emph{descending} if it contains consecutive rows $i,i+1$ such that all 1s in row $P(i,\cdot)$ strictly precede all 1s in row $P(i+1,\cdot)$.  It is \emph{ascending} if there exists $i,i+1$ such that all 1s in row $P(i+1,\cdot)$ strictly precede all 1s in row $P(i,\cdot)$.
A pattern $P$ is \emph{legal} if it
is not both descending and ascending.
\end{definition}

Pettie's grafting theorem~\cite{Pettie-SoCG11} implies that if $P_0$ is legal, 
$\ex(n,P_0) \leq n\cdot f_0(n)$,
$\ex(n,P_1) \leq n\cdot f_1(n)$,
and $f_0(n),f_1(n) < \mbox{polylog}(n)$,
then 
$\ex(n,P_0\otimes P_1)=O(n\cdot f_0(n)f_1(n))$.\footnote{The statement proved in~\cite{Pettie-SoCG11} only shows that if $P_0,P_1$ are linear, then $P_0\otimes P_1$ is too, but it easily extends to non-linear patterns.}

\subsection{New results}

In this paper we give three new pattern transformations,
which yield new infinite classes of linear patterns.
These results support \Cref{conj:PT}
and have implications for the classification of
weight-5 patterns.

\restateable{theorem}{threeones}{\label{thm:three-ones}
Let $P\in \{0,1\}^{k\times k}$ 
be any permutation pattern. 
Define $P' \in \{0,1\}^{(k+2)\times (k+2)}$
to be identical to $P$ on indices 
$([1,w]\cup [w+3,k+2]) \times [2,k+1]$, 
for some $w\in\{0,\ldots,k\}$,
with additional 1s at $P'(w+1,1)=P'(w+2,1)=P'(w+2,k+2)=1$,
and all other entries 0.  
(See the diagram below.)
Then $\ex(n,P')=O(n)$.
\[
P' = \left(\begin{array}{ccccc} \cline{2-4}
 & \multicolumn{1}{|c}{} && \multicolumn{1}{c|}{} &\\
 & \multicolumn{1}{|c}{} & \rb{2}{$P$} & \multicolumn{1}{c|}{}\\\cline{2-4}
 \bu\\
 \bu &&&& \bu\\
 \cline{2-4}
 & \multicolumn{1}{|c}{} && \multicolumn{1}{c|}{} &\\
 & \multicolumn{1}{|c}{} & & \multicolumn{1}{c|}{}\\\cline{2-4}
\end{array}\right).
\]
}

\restateable{theorem}{twoextension}{\label{thm:two-sided-extension}
Suppose $P\in\{0,1\}^{k\times l}$ is legal and 
there exists a row $w$ with $P(w,1)=P(w,l)=1$.  Obtain $P' \in \{0,1\}^{(k+1)\times (l+2)}$ by putting $P$ in the submatrix
on indices $([1,w-1]\cup [w+1,k+1])\times [2,l+1]$
and inserting a weight-2 row with $P'(w,1)=P'(w,l+2)=1$.
(See the diagram below.)
Then
\[
 \ex(n,m,P')\le \ex(3n,4m,P) + 2n + 4m.
\]
}

\[
P' = \left(\begin{array}{cccccc} \cline{2-5}
 & \multicolumn{1}{|c}{} &&& \multicolumn{1}{c|}{} &\\
 & \multicolumn{1}{|c}{} && \zero{\hcm[-.32]\rb{2}{$P$}} & \multicolumn{1}{c|}{}\\\cline{2-5}
\bu &&&&& \bu\\\cline{2-5}
 & \multicolumn{1}{|c}{\bu} &&& \multicolumn{1}{c|}{\bu} &\\
  & \multicolumn{1}{|c}{} &&& \multicolumn{1}{c|}{} &\\\cline{2-5}
\end{array}\right).
\]

Naturally \Cref{thm:two-sided-extension} can be applied to insert a row after $w$, by considering the reflection of $P$.
There is also a one-sided version of \Cref{thm:two-sided-extension} that only adds a single 1 but applies to a wider class of patterns $P$.

\restateable{theorem}{oneextension}{\label{thm:one-sided-extension}
Suppose $P\in\{0,1\}^{k\times l}$ is legal and there 
exists a row $w$ with $P(w,l)=1$.  Obtain $P' \in \{0,1\}^{(k+1)\times (l+1)}$ by putting $P$ in the 
submatrix on indices $([1,w-1]\cup [w+1,k+1])\times [1,l]$
and inserting a weight-1 row with $P'(w,l+1)=1$.  (See the diagram below.)
Then
\[
 \ex(n,m,P')\le \ex(2n,3m,P) + n+m.
\]
\[
P' = \left(\:\begin{array}{ccccc} \cline{1-4}
 \multicolumn{1}{|c}{} &&& \multicolumn{1}{c|}{} &\\
 \multicolumn{1}{|c}{} && \zero{\hcm[-.32]\rb{2}{$P$}} & \multicolumn{1}{c|}{}\\\cline{1-4}
 &&&& \bu\\\cline{1-4}
 \multicolumn{1}{|c}{} &&& \multicolumn{1}{c|}{\bu} &\\
  \multicolumn{1}{|c}{} &&& \multicolumn{1}{c|}{} &\\\cline{1-4}
\end{array}\right).
\]
}

\Cref{thm:two-sided-extension} and 
\Cref{thm:one-sided-extension} are very powerful.
They can reprove the linearity of all known linear patterns except for (double) permutations~\cite{MarcusT04,Geneson09}
and $L_P$~\cite{Pettie-GenDS11}.
\Cref{thm:two-sided-extension} directly 
gives yet another proof
that monotone patterns~\cite{Keszegh09,Pettie-SoCG11} 
are linear and that Fulek's pattern $L_F$ is linear.
\[
\patbits{010/111} 
\xrightarrow{\,\mbox{\small Thm.~\ref{thm:two-sided-extension}}\,}
\patbits{00100/10001/01110}
\]
We can apply 
\Cref{thm:two-sided-extension} 
twice, once to add a row and another to add a column, to yield a 
\emph{new} linear pattern 
containing Tardos's $L_T$.
\[
\patbits{10/11}
\xrightarrow{\,\mbox{\small Thm.~\ref{thm:two-sided-extension}}\,}
\patbits{0100/1001/0110}
\xrightarrow{\,\mbox{\small Thm.~\ref{thm:two-sided-extension}}\,}
\patbits{00100/01000/10001/01010/00100}
\]
Similarly, the linearity of 
Geneson's pattern~\cite{Geneson26}
follows from \Cref{thm:one-sided-extension}.
\[
 \patword{12211}
 \xrightarrow{\,\mbox{\small Thm.~\ref{thm:one-sided-extension}}\,}
 \patword{133112}.
\]
\Cref{thm:one-sided-extension} can also be used to classify previously unclassified \emph{non}linear patterns.  One such example is
\[
\patword{1212}
 \xrightarrow{\,\mbox{\small Thm.~\ref{thm:one-sided-extension}}\,}
\patword{13132}
\]
Therefore the extremal function of the last pattern is asymptotically identical to the first, namely $\Theta(n\alpha(n))$.

\medskip
Pettie and Tardos~\cite{PettieT25} identified several classes of $\mathcal{Q}$-free patterns that must all be linear if \Cref{conj:PT} is true.  
Repeated application of \Cref{thm:two-sided-extension} implies that the $\{O_t\}$ patterns of~\cite{PettieT25} are all 
linear, and in fact all the doublings of the $\{O_t\}$ patterns.  

\[
\begin{aligned}
 O_2&=\patwordc{2112},& O_3&=\patwordc{213312},\\
 O_4&=\patwordc{32411423},& O_5&=\patwordc{3241551423}.
\end{aligned}
\]

\restateable{corollary}{oscillating}{\label{cor:Ok-family}
For every $t\geq 2$, 
the oscillating pattern $O_t$ has $\ex(n,O_t) = O(n)$.
}

\Furedi{} and Hajnal~\cite{FurediH92}
proved that $\ex(n,O_2)=O(n)$, but no other oscillating patterns were previously classified.

\subsubsection{The Classification of Weight-5 Patterns}

Theorems \ref{thm:three-ones}--\ref{thm:one-sided-extension}
give us a complete classification of
linear light weight-5 patterns, and in fact all linear weight-5 patterns. 

\begin{corollary}\label{cor:light-weight-five}
    The extremal function of every light 
    weight-5 pattern is either $O(n)$ or $\Omega(n\alpha(n))$.  Such a $P$ is linear if and only if it is $\mathcal{Q}$-free.
\end{corollary}

Before this work, the only light weight-5 patterns whose linear/nonlinear status was unknown were the three below.
\begin{align*}
Q_{41423} &= \patword{41423}
&&& Q_{41432} &= \patword{41432}
&&& Q_{42314} &= \patword{42314}
\end{align*}

One can actually deduce that $Q_{41423}$ and $Q_{41432}$ are linear by grafting a rotation of $O_2$ onto the side of a subpattern of $L_T$
and applying trivial transformations.
\[
\patbits{010/001/101}
 \xrightarrow{\,\mbox{\small graft}\,}
\patbits{0100/0010/0001/0001/1010}.
\]
\Cref{thm:three-ones} implies 
that $Q_{42314}$ is linear.

\medskip 

We are also able to classify all weight-5 patterns according to their linear/nonlinear status.

\begin{corollary}\label{cor:weight-five}
    If $P$ is a weight-5 nonlinear pattern,
    its nonlinearity can be deduced from the nonlinearity of $Q_3'$ or $R_4$ via trivial transformations.
\begin{align*}
Q_3' &= \patword{1323},
&&& R_4 &= \patbits{1000/0001/0001/0110}.
\end{align*}
\end{corollary}

After resolving the linearity of all light patterns in \Cref{cor:light-weight-five}, 
the only \emph{non}light weight-5 patterns whose linear/nonlinear status 
was outstanding were
$L_1$ and $L_2$.  \Cref{thm:three-ones} implies that both are linear.
\begin{align*}
    L_{1} &= \patbits{0010/0100/1000/1001},
&&&    L_2 &= \patbits{0100/0010/1000/1001}.
\intertext{Aside from the new patterns shown to be linear, 
we also prove sharp $\Theta(n\alpha(n))$ 
bounds on two of the outstanding weight-5 patterns
$\hat{Q}_3,\tilde{Q}_3$.}
    \hat{Q}_3 &= \patbits{1010/0100/0101},
&&&   \tilde{Q}_3 &= \patword{31231}.
\end{align*}

\begin{restatable}{theorem}{Qhat}\label{thm:Qhat-alpha}
    $\ex(n,m,\hat{Q}_3) = \Theta(n\alpha(n,m)+m)$.
\end{restatable}

\begin{restatable}{theorem}{Qtilde}\label{thm:Qtilde-alpha}
    $\ex(n,m,\tilde{Q}_3) = \Theta(n\alpha(n,m)+m)$.
\end{restatable}

\subsubsection{Minimally Nonlinear Patterns}

Keszegh~\cite{Keszegh09} proved that there are infinitely many minimally nonquasilinear patterns\footnote{A function $f$ is \emph{quasilinear} if it is $n2^{O(\alpha^t(n))}$ for some $t=O(1)$.} 
with respect to $\prec$ containment.  
He was unable to extend this statement to 
minimally \emph{nonlinear} patterns because
the status of double permutation patterns 
was still open at the time.  
Geneson~\cite{Geneson09} 
proved that double permutations are linear, 
implying the existence of infinitely many minimally nonlinear patterns.  See Pettie~\cite{Pettie-FH11} for a simplified presentation of Keszegh's proof.
Further structural bounds on minimally nonlinear matrices, 
such as linear bounds on their weight 
with respect to their dimensions, 
were obtained by CrowdMath~\cite{CrowdMath2018} 
and Geneson and Tsai~\cite{GenesonTsai20}.

The prior results~\cite{Keszegh09,Geneson09,Pettie-FH11,CrowdMath2018,GenesonTsai20} 
still do not demonstrate that a 
\emph{specific} infinite class of patterns 
is minimally nonlinear.  Keszegh~\cite{Keszegh09} conjectured that every pattern in 
the family $(H_k)_{k\geq 0}$ is minimally nonlinear.
We define $H_k$ formally in \Cref{sec:hk-minimal}
and list the first three patterns in the series below.
Keszegh~\cite{Keszegh09} proved that 
$H_0$ 
is minimally nonlinear and that
all the patterns $(H_k)_{k\geq 0}$ 
have extremal function $\ex(n,H_k)=\Theta(n\log n)$,
once the linearity of 
double permutations was resolved~\cite{Geneson09}.

\[
H_0 = \scalebox{.9}{$\patbits{0110/0001/0001/1000}$}
\quad
H_1 = \scalebox{.9}{$\patbits{0110000/0000010/0000100/1000000/0000001/0000001/0001000}$},
\quad
H_2 = \scalebox{.9}{$\patbits{0110000000/0000010000/0000100000/1000000000/0000000010/0000000100/0001000000/0000000001/0000000001/0000001000}$}.
\]

We prove Keszegh's conjecture that $(H_k)_{k\geq 0}$ is minimally nonlinear, and also that 
all minimally nonlinear patterns $P$, with the exception of $C_4$, have $\ex(n,P)=O(n\log n)$.

\restateable{theorem}{hkminimal}
{\label{thm:hk-minimal}
For every $k\geq 0$, the pattern $H_k$ is minimally nonlinear.
}

\restateable{theorem}{minimallynonlinear}{\label{thm:minimally-nonlinear}
    If $P\neq C_4$ is minimally nonlinear, then $\ex(n,P)=O(n \log n)$.
}

\subsubsection{Extremal Function $\Theta(n\log n/\log\log n)$}

We provide \emph{general} methods for obtaining 
$O(n\log n/\log\log n)$ upper bounds and $\Omega(n\log n/\log\log n)$ lower bounds for any pattern with certain structural properties.
Ten weight-6 patterns satisfy both the upper and lower bounds.
We give two examples below.
\[
\patbits{00110/01001/00100/10000},
\qquad\qquad
\patbits{01010/00100/10101}.
\]

\subsection{Organization}

\Cref{sect:linear-patterns} proves Theorems~\ref{thm:three-ones}--\ref{thm:one-sided-extension}, which give us infinitely many more linear patterns.
\Cref{sect:minimally-nonlinear-patterns} contains all the new results on minimally nonlinear patterns, namely an explicit family of infinitely many minimally nonlinear patterns, 
exponential lower bounds on the number of such patterns in terms of their weight, and an $O(n\log n)$ upper bound on the extremal function of such patterns.
\Cref{sect:logn-over-loglogn} gives two general methods for proving $\Omega(n\log n/\log\log n)$ lower bounds and $O(n\log n/\log\log n)$ upper bounds on patterns satisfying a few structural properties.  Ten weight-6 patterns are identified with sharp $\Theta(n\log n/\log\log n)$ bounds.
Appendix~\ref{sect:new-quasilinear-bounds} 
proves \Cref{thm:Qhat-alpha,thm:Qtilde-alpha},
giving $\Theta(n\alpha(n,m)+m)$ bounds on 
$\hat{Q}_3$-free and $\tilde{Q}_3$-free matrices.

\medskip 

We conclude in \Cref{sec:conclusion} with a survey of the remaining small unclassified patterns.

\section{New Pattern Transformations}\label{sect:linear-patterns}

We prove \Cref{thm:three-ones} in \Cref{sec:theorem-threeones}
and \Cref{thm:two-sided-extension,thm:one-sided-extension} in \Cref{sec:reductions}.
These theorems yield infinitely more linear patterns, and sharp bounds on some nonlinear patterns.

\subsection{Proof of \Cref{thm:three-ones}}\label{sec:theorem-threeones}

The proof uses a few facts summarized in \Cref{thm:FH-SW-Ge}.

\begin{theorem}[See \cite{MarcusT04,Geneson09,Cibulka09,CibulkaK17,Fox13,Klazar00}]\label{thm:FH-SW-Ge}
For any $k$-permutation $\pi$ let $P_\pi$ be its 0--1 matrix and $P_\pi^{\dbl}$ be obtained from $P_\pi$ by doubling every column.  There are constants 
$C^{\FH}_\pi, C^{\SW}_\pi, C^{\Ge}_\pi$ 
such that the following hold.
\begin{enumerate}
    \item $\ex(n,m,P_\pi) \leq C^{\FH}_\pi(n+m)$.\label{item:FH}
    \item $\ex(n,m,P_\pi^{\dbl}) \leq C^{\Ge}_\pi(n+m)$.\label{item:Ge}
    \item The number of $m$-permutations that avoid $\pi$ is at most $(C^{\SW}_\pi)^m$.
\end{enumerate}
\end{theorem}

The constants $C^{\FH}_\pi, C^{\SW}_\pi, C^{\Ge}_\pi$ are known to be $2^{O(k)}$, and $2^{\Omega(\sqrt{k})}$ for some permutation $\pi$~\cite{Fox13,Cibulka09,GenesonT17},
and $2^{\tilde{O}(\sqrt{k})}$ for nearly all permutations $\pi$~\cite{Geneson26b,GenesonT17}.  
We use $C^{\FH}_k,C^{\SW}_k,C^{\Ge}_k$ to be the maximum
over all $k$-permutations $\pi$.

\paragraph{The Patterns.}
Following the terminology of \Cref{thm:three-ones},
let $\pi$ be a $k$-permutation, $P=P_\pi$ be its 0--1 matrix,
and $P'$ be the target matrix, for which we will prove that $\ex(n,P')=O(n)$.
Define $P''$ to be $P'$ after removing the weight-1 row outside of $P$.  
Obtain $L$ and $R$ by removing the last and first column of $P''$, respectively

\[
P' = \scalebox{.85}{$\displaystyle \left(\begin{array}{ccccc} \cline{2-4}
 & \multicolumn{1}{|c}{} && \multicolumn{1}{c|}{} &\\
 & \multicolumn{1}{|c}{} & \rb{2}{$P$} & \multicolumn{1}{c|}{}\\\cline{2-4}
 \bu\\
 \bu &&&& \bu\\
 \cline{2-4}
 & \multicolumn{1}{|c}{} && \multicolumn{1}{c|}{} &\\
 & \multicolumn{1}{|c}{} & & \multicolumn{1}{c|}{}\\\cline{2-4}
\end{array}\right)$},\quad
P'' = \scalebox{.85}{$\displaystyle\left(\begin{array}{ccccc} \cline{2-4}
 & \multicolumn{1}{|c}{} && \multicolumn{1}{c|}{} &\\
 & \multicolumn{1}{|c}{} & \rb{2}{$P$} & \multicolumn{1}{c|}{}\\\cline{2-4}
 \bu &&&& \bu\\
 \cline{2-4}
 & \multicolumn{1}{|c}{} && \multicolumn{1}{c|}{} &\\
 & \multicolumn{1}{|c}{} & & \multicolumn{1}{c|}{}\\\cline{2-4}
\end{array}\right)$},\quad
L = \scalebox{.85}{$\displaystyle\left(\begin{array}{ccccc} \cline{2-4}
 & \multicolumn{1}{|c}{} && \multicolumn{1}{c|}{} &\\
 & \multicolumn{1}{|c}{} & \rb{2}{$P$} & \multicolumn{1}{c|}{}\\\cline{2-4}
 \bu &&&& \\
 \cline{2-4}
 & \multicolumn{1}{|c}{} && \multicolumn{1}{c|}{} &\\
 & \multicolumn{1}{|c}{} & & \multicolumn{1}{c|}{}\\\cline{2-4}
\end{array}\right)$},\quad
R = \scalebox{.85}{$\displaystyle\left(\begin{array}{ccccc} \cline{2-4}
 & \multicolumn{1}{|c}{} && \multicolumn{1}{c|}{} &\\
 & \multicolumn{1}{|c}{} & \rb{2}{$P$} & \multicolumn{1}{c|}{}\\\cline{2-4}
 &&&& \bu\\
 \cline{2-4}
 & \multicolumn{1}{|c}{} && \multicolumn{1}{c|}{} &\\
 & \multicolumn{1}{|c}{} & & \multicolumn{1}{c|}{}\\\cline{2-4}
\end{array}\right)$}.
\]

Before proving \Cref{thm:three-ones} we 
introduce two useful lemmas.

\begin{lemma}\label{lem:double-perm}
Let $s \in [N]^{N(N+1)}$ be a string in 
which each symbol $i\in [N]$ appears 
at least $N+1$ times.  Then $s$ contains
a subsequence $s' = \pi(1)\pi(1)\pi(2)\pi(2)\cdots \pi(N)\pi(N)$ for some permutation $\pi : [N]\to[N]$.
\end{lemma}

\begin{proof}
    When $N\geq 1$ every symbol 
    appears at least twice.
    Let $i$ be the first repeated symbol in $s$, i.e.,
    $s = \alpha \, i\,  \beta\,  i\, \cdots$, where no symbol appears twice in $\alpha\beta$.
    Set $\pi(1) = i$ and remove
    $\alpha i\beta i$ and all subsequent 
    occurrences of $i$.  The remaining string satisfies the preconditions of the lemma 
    with an alphabet size of $N-1$,
    so the claim follows by induction.
    The base case $N=0$ is vacuous.
\end{proof}

For any $N$-permutation $\sigma$, 
define $S_\sigma$ to be the $N\times 2N$ matrix 
$(P_\sigma \,|\, I)$, where $I$ is the identity.

\begin{lemma}\label{lem:cap-forcing}
There exists an $N=N(k)$ and an $N$-permutation $\sigma^*$
such that $P''$ is contained in \emph{every} 
pattern obtained by shuffling the rows of $S_{\sigma^*}$.
Moreover, $P_{\sigma^*}$ contains a 1 in the top right corner.
\end{lemma}

\begin{proof}
Suppose a row-shuffle $(X\,|\,Y)$ of $S_\sigma$ is $P''$-free.  
It must be that its left half $X$ is $L$-free and that its right half $Y$ is $R$-free.
By \Cref{thm:FH-SW-Ge}(3) there are at most
$(C_{k+1}^{\SW})^N (C_{k+1}^{\SW})^N$ pairs 
$(X\,|\, Y)$ having this property.
Observe that each $(X\,|\, Y)$ uniquely determines $\sigma$,
since $Y$ was initially $I$ before shuffling.
Thus, there are at most $(C_{k+1}^{\SW})^{2N}$ 
permutations $\sigma$
for which $S_\sigma$ admits some $P''$-free row-shuffling.
We choose $N$ such that $(N-1)! > (C_{k+1}^{\SW})^{2N}$.
By the pigeonhole principle, there is an $(N-1)$-permutation
$\hat{\sigma}$ such that every row-shuffling of 
$S_{\hat{\sigma}}$ contains $P''$.  We let $\sigma^*$
be an $N$-permutation such that $P_{\sigma^*}$ is $P_{\hat{\sigma}}$ with a new 1 placed in the NE corner.
\end{proof}

We are now prepared to prove \Cref{thm:three-ones}.

\threeones*

\begin{proof}
Let $A$ be an $n\times m$, $P'$-free matrix, and let $A(r,\last(r))=1$ be the last 1 in row $r$.
Define a new $m\times m$ matrix $B$ such that
\[
B(z,j) = 1 
\quad\Longleftrightarrow\quad
 j < z \mbox{ and for some row $r$, $\last(r)=z$ and $A(r,j)=1$.}
\]

Let $n_z$ be the number of rows $r$ with $\last(r)=z$
and let $A_z$ be the submatrix of $A$ 
on those $n_z$ rows, with column $z$ zeroed out.
Let $m_z$ be the number of columns that are 
non-zero in $A_z$.
It follows that $\sum_z n_z \leq n$ and 
$\sum_z m_z = |B|$.

Obtain $\tilde{L}$ by repeating the 
row of $L$ outside of $P$,
or equivalently, removing the last 
column of $P'$.
Observe that every $A_z$ is $\tilde{L}$-free, 
for any such occurrence extends to an occurrence of 
$P'$ in $A$ by including a 1 from column $z$.
It follows from \Cref{thm:FH-SW-Ge}(\ref{item:Ge}) that 
\begin{align*}
    |A_z| &\leq m_z + C^{\Ge}_{k+1}(n_z+m_z),
\intertext{and hence}
|A| &= \sum_z (n_z + |A_z|) \leq 
(C^{\Ge}_{k+1} +1)(n + |B|).
\end{align*}
It remains to bound $|B|$.  Fix $N,\sigma^*$ from \Cref{lem:cap-forcing} and let $P_{\sigma^*}$ be the permutation matrix of $\sigma^*$.
Let $P_{\sigma^*}^{\times}$ be the 
$N(N+1)\times N$ matrix 
obtained by repeating every row of $P_{\sigma^*}$ $N+1$ times.
We claim that $B$ is 
$P_{\sigma^*}^{\times}$-free, 
hence
\[
|B| 
\leq \ex(m,m,P_{\sigma^*}^{\times}) 
\leq N\cdot \ex(m,m,P_{(\sigma^*)^{-1}}^{\dbl})
\leq N\cdot C^{\Ge}_{N}(m+m).
\]
Suppose that $P_{\sigma^*}^{\times}$ appears in $B$
on rows $z_1 < \cdots < z_{N(N+1)}$ and 
columns $c_1<\cdots <c_N$.  Since the NE corner of $P_{\sigma^*}^{\times}$ is 1, $B(z_1,c_{N})=1$ implies 
$c_N < z_1$,
so in $A$, the column order is 
$c_1 < \cdots < c_N < z_1 < \cdots <z_{N(N+1)}$.
Let $r_1 < r_2 < r_3 < \cdots < r_{N(N+1)}$
be the rows in $A$ witnessing the 1s of $P_{\sigma^*}^{\times}$ in $B$.
Note that there is \emph{no correspondence}
between the order of $z_1,\ldots,z_{N(N+1)}$
and the order of $r_1,\ldots,r_{N(N+1)}$.
By \Cref{lem:double-perm}, 
$(r_i)_{i\in [N(N+1)]}$
contains a subsequence $(r'_i)_{i\in [2N]}$
where each row of $P_{\sigma^*}$
is identified with two consecutive
rows $r'_{2i-1},r'_{2i}$.
In other words, $A$ contains 
some row-shuffling of $S_{\sigma^*}$.
The diagram below illustrates the situation where the witness rows $(r'_i)$ in $A$ 
\emph{happen to coincide} with the 
order of their counterparts in $(z_i)$.
For clarity we reindex these counterparts 
as $z_1 < \cdots < z_{2N}$.
\begin{align*}
    \text{In $B$:} 
    &\begin{array}{r}
    z_1\\
    z_2\\
    \ \\
    \ \\
    \zero{\vdots}\;\:\\
    \ \\
    \ \\
    z_{2N}\\
    \end{array}
    \left(\begin{array}{c@{\quad}c@{\quad}c@{\quad}c}
    \zero{\rb{4.5}{\hcm[-.05]$c_1$}}
    &
    \zero{\rb{4.5}{\hcm[-.05]$c_2$}}&
    \zero{\rb{4.5}{\hcm[-.14]$\cdots$}}
    &
    \zero{\rb{4.5}{\hcm[0]$c_N$}}\bu\\
    &&&\bu\\
    &\bu&&\\
    &\bu&&\\
    &&\bu&\\
    &&\bu&\\
    \bu&&&\\
    \bu&&&
    \end{array}\right)
    &&& \text{In $A$:} &
    \begin{array}{r}
    \ \\
    r'_1\\
    r'_2\\
    \ \\
    \ \\
    \zero{\vdots}\;\:\\
    \ \\
    \ \\
    r'_{2N}\\
    \end{array}
    \left(\begin{array}{cccccccccccc}\:
    \zero{\rb{4.5}{$c_1$}} &
    \zero{\rb{4.5}{$c_2$}} &
    \zero{\rb{4.5}{\hcm[-.05]$\cdots$}} &
    \zero{\rb{4.5}{\hcm[0.1]$c_N$}}\bu &
    \zero{\rb{4.5}{\hcm[0.12]$z_1$}}\bu&
    \zero{\rb{4.5}{\hcm[0.08]$z_2$}} &
    &
    &
    \zero{\rb{4.5}{\hcm[.35]$\cdots$}}&
    &
    &\istrut{10}
    \zero{\rb{4.5}{\hcm[-.3]$z_{2N}$}}\\
    &&&\bu&&\bu&&&&&&\\
    &\bu&&&&&\bu&&&&&\\
    &\bu&&&&&&\bu&&&&\\
    &&\bu&&&&&&\bu&&&\\
    &&\bu&&&&&&&\bu&&\\
    \bu&&&&&&&&&&\bu&\\
    \bu&&&&&&&&&&&\bu
    \end{array}\:\right)
\end{align*}
The $i$th pair of 
duplicate rows 
of $P_{\sigma^*}^{\times}$
is mapped to two rows 
$u_{\lo}(i) < u_{\hi}(i)$ 
in $A$, where both are among the 
$\{r'_i\}_{i\in [2N]}$.
Therefore
rows $u_{\lo}(1),\ldots,u_{\lo}(N)$
of $A$
contain a row-shuffling
of $S_{\sigma^*}$, 
and according to 
\Cref{lem:cap-forcing}, a copy of $P''$. 
Let $i$ be such that
row $u_{\lo}(i)$ contains the weight-2
row of $P''$.
If we insert row $u_{\hi}(i) > u_{\lo}(i)$,
this puts a 1 immediately below the left 1
in row $u_{\lo}(i)$, and another 1
immediately before or after the right 1 
in row $u_{\lo}(i)$, the two possibilities 
being marked by $\circ$ in the diagram below.
This creates a copy of $P'$ in $A$,
contradicting its $P'$-freeness, 
which means that $B$ is 
$P_{\sigma^*}^{\times}$-free.

\[
P' \prec \;
\begin{array}{r}
\ \\
\ \\
u_{\lo}(i)\\
u_{\hi}(i)\istrut{6}\\
\ \\
\ \\
\end{array}
\left(\begin{array}{ccccccc} \cline{2-4}
 & \multicolumn{1}{|c}{} && \multicolumn{1}{c|}{} &\\
  & \multicolumn{1}{|c}{} & \rb{2}{$P$} & \multicolumn{1}{c|}{}\\\cline{2-4}
 \bu &&&&& \bu &\\
 \bu\istrut{6} &&&& \circ\zero{\rb{4.5}{\hcm[.35]$\swarrow$}} && \circ\zero{\rb{4.5}{\hcm[-.75]$\searrow$}}\\
\cline{2-4}
 & \multicolumn{1}{|c}{} && \multicolumn{1}{c|}{} &\\
  & \multicolumn{1}{|c}{} &  & \multicolumn{1}{c|}{}\\\cline{2-4}
\end{array}\right).
\]

The known upper bounds on $C_{k+1}^{\FH}$ and $C_{k+1}^{\SW}$~\cite{Fox13,Cibulka09} 
imply that $N(k) = 2^{\Theta(k)}$
suffices for \Cref{lem:cap-forcing}.
Since $C_{N}^{\Ge} = 2^{O(N)}~\cite{GenesonT17}$,
we can then bound $|A|$ by:
\begin{align*}
    |A| &\leq (C_{k+1}^{\Ge}+1)(n+N\cdot C_N^{\Ge}(m+m)) 
    = 2^{O(k)}n + 2^{2^{O(k)}}m.
\end{align*}
\end{proof}

\subsection{Proofs of \Cref{thm:two-sided-extension,thm:one-sided-extension}}\label{sec:reductions}

The proofs in this section build on the 
\emph{visibility graph} technique of Fulek~\cite{Fulek09}
and Geneson and Shen~\cite{GenesonShen}.
The transformations of \Cref{thm:two-sided-extension,thm:one-sided-extension}
only apply to \emph{legal} patterns $P$ 
(see \Cref{def:legal}), meaning that $P$ is not 
both ascending and descending.  
Because \Cref{thm:one-sided-extension}
lacks the symmetry of \Cref{thm:two-sided-extension},
the proof is somewhat more complex.  We begin with it.

\oneextension*

\begin{proof}
If $w=1$ then $\ex(n,m,P')\leq \ex(n,m,P)+n+m$,
so we shall assume $w>1$.
Let $A_0$ be the $P'$-free $n\times m$ matrix,
and $A$ be derived by removing the last 1 from each row of $A_0$.
Define $\last(r)$ to be the column of the last 1 in row $r$ of $A$.
Attribute to row $r$ the \emph{bar} 
$(-\infty,\last(r)+1/2]$.
We say two bars are \emph{visible} if one can draw a vertical strut between them that intersects only those 
bars.  The \emph{visibility graph} 
identifies bars with vertices and visible pairs with 
edges.  It is planar~\cite{DuchetHVM83,TamassiaT86,Fulek09}, 
and in our case 3-colorable.  
Fix any 3-coloring of the rows.

At any $x$-coordinate less than $\last(r)+1/2$, 
row $r$ has a unique upper neighbor.
Partition each row into \emph{pieces} for which the 
identity of the upper neighbor is constant, 
with a sentinel $(-\infty,\infty)$ at row 0.
Each piece $u$ is identified with its visible pair $(\Top(u),\Bottom(u))$, where $\Bottom(u)$ is the 
row of the piece in $A$ 
and $\Top(u)$ is its constant upper neighbor.
The piece $u$ retains the color of $\Bottom(u)$.
Derive $A'$ from $A$ by putting each piece in a 
distinct row, without changing the columns.  
The order of the pieces depends on whether 
$P$ is nonascending or nondescending.
\begin{description}
    \item[Nondescending $P$.] 
    The pieces are ordered lexicographically by $(\Bottom(u), -\Top(u))$.
    \item[Nonascending $P$.] 
    The pieces are ordered lexicographically by $(\Top(u), -\Bottom(u))$.
\end{description}

\begin{figure}
    \centering
\begin{tabular}{c}
    \includegraphics[width=0.98\linewidth]{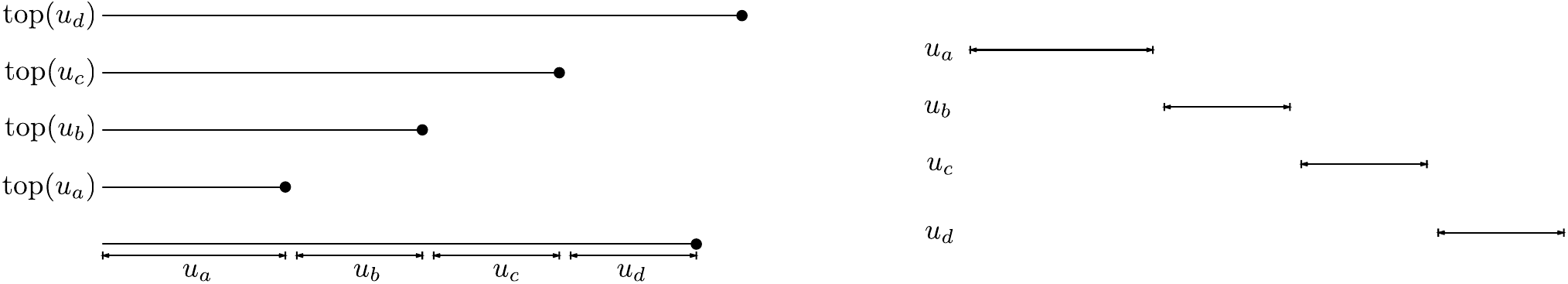}\\
    (A)\\\ \\
    \includegraphics[width=0.98\linewidth]{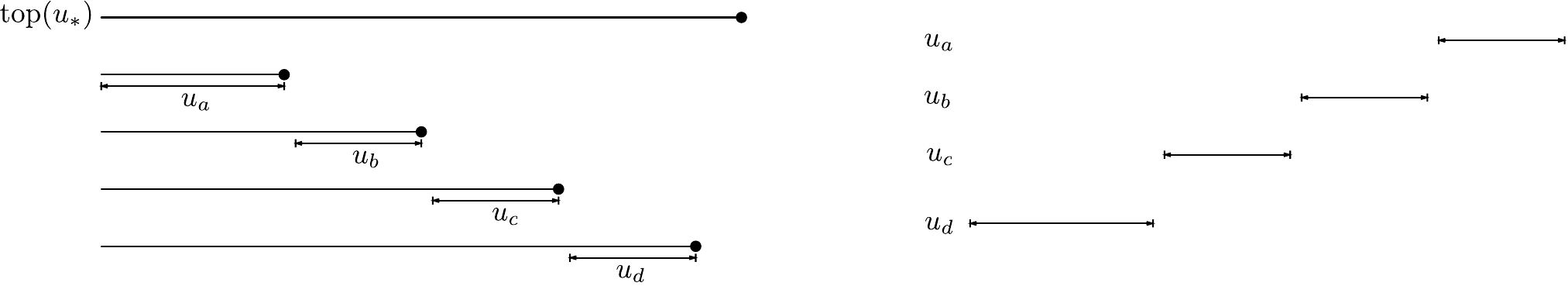}\\
    (B)
\end{tabular}
    \caption{(A) Left: pieces $u_a,\ldots,u_d$ with a common $\Bottom$-row.  
    Right: when $P$ is nondescending, 
    the pieces of a common row in $A$ are mapped to distinct rows in $A'$, in descending order. (B) Left: pieces $u_a,\ldots,u_d$ with a common $\Top$-row.  
    Right: when $P$ is nonascending, the pieces $u_a,\ldots,u_d$ in $A$ are mapped to distinct rows in $A'$, in increasing order.}
    \label{fig:pieces}
\end{figure}

In either case, let $A_c$ be the $n_c\times m$ 
submatrix of $A'$ on the rows colored $c$.  
A row has one piece starting at the first column
and every subsequent piece can be charged to the termination of some bar.
It follows that $n_1+n_2+n_3 \leq 2n$.
We now claim that each $A_c$ is $P$-free.

\paragraph{Case: Nondescending $P$.}  
A row of $A$ colored $c$ is represented by a 
\emph{contiguous} interval of pieces in $A_c$.  
The tie-breaking value for pieces with a common $\Bottom(u)$
is $-\Top(u)$, which places them in descending order.  See \Cref{fig:pieces}(A).
Since $P$ is nondescending, 
it follows that any occurrence of $P$ in $A_c$ is 
mapped to \emph{distinct} rows in $A$.  
Let $u_{w-1},u_w$ be the pieces holding rows $w-1,w$ in an occurrence of $P$ in $A_c$, 
and $x < y$ be their original rows in $A$.
Let $j$ be the column of the first 1 in $u_w$.
$P$ is nondescending, so the last 1 in row $u_{w-1}$
must be at or to the right of $j$.  It cannot be that
$x$ is visible from $y$ since they have the same color,
so $x < \Top(u_w) < y$.  
Recall that $\Top(u_w)$ is the upper neighbor of $y$ throughout $u_w$.  
The last 1 in row $\Top(u_w)$ of $A_0$ extends the occurrence of $P$ in $A_c$ to an occurrence of $P'$ in $A_0$, a contradiction.  Thus, $A_c$ is $P$-free.

\paragraph{Case: Nonascending $P$.} 
In this case the pieces originating from a color-$c$
row in $A$ are \emph{not} necessarily assigned to a contiguous interval of rows in $A_c$
since they are sorted primarily by 
their $\Top$ coordinate, not $\Bottom$.
Nonetheless, we argue that an occurrence of $P$ in $A_c$ implies an occurrence of $P'$ in $A_0$.

First observe that the pieces 
originating from the same color-$c$
row in $A$ are arranged in ascending order in $A_c$
due to ordering first by $\Top$.
See \Cref{fig:pieces}.
Now fix an occurrence of $P$ in $A_c$
and consider adjacent rows $j-1,j$ in $P$, 
their corresponding pieces $u_{j-1} < u_j$ in $A_c$,
and their original rows $x,y$ in $A$.  
(We also use $u_j$ to refer to the index of the row in $A'$ holding $u_j$.)
Suppose, for the purpose of obtaining a contradiction, 
that the order in $A$ does not respect the order in $A_c$,
i.e., $x\geq y$.  It cannot be that $x=y$ 
since pieces from the same row are in ascending order 
and $P$ is nonascending.
If $y<x$, then $y$'s bar $(-\infty,\last(y))$ 
at least partially occludes $u_{j-1}$,
but since $x,y$ have the same color $c$, 
they cannot be adjacent.  It follows
that $\Top(u_j) < \Top(u_{j-1})$ and hence $u_j < u_{j-1}$, a contradiction.  The conclusion here is 
only that an occurrence of $P$ in $A_c$ 
corresponds to an occurrence of $P$ in $A$.

\begin{figure}
    \centering
    \begin{tabular}{c@{\hcm}c}
    \includegraphics[width=0.46\linewidth]{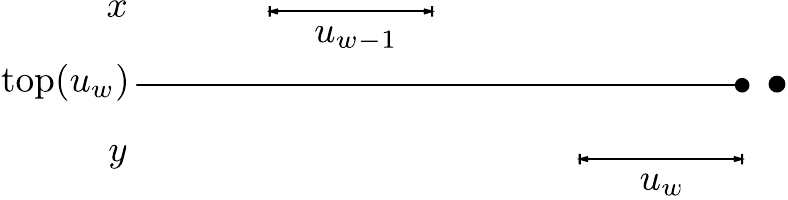} &
    \includegraphics[width=0.46\linewidth]{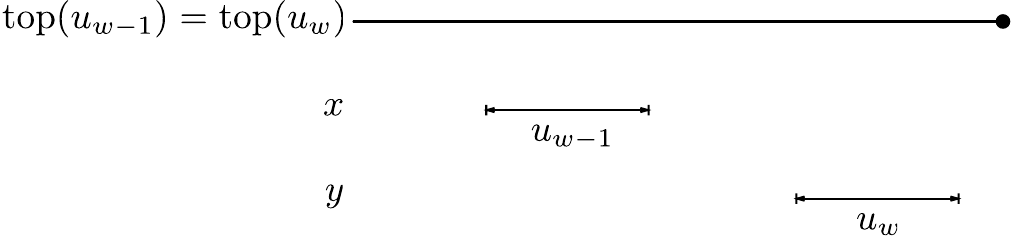}\\
    (A) & (B)
    \end{tabular}
    \caption{Left: When $\Top(u_w)$ lies between rows $x$ and $y$, an occurrence of $P$ in $A_c$ extends to $P'$ in $A_0$ with the last 1 of row $\Top(u_w)$.
    Right: It is impossible for $\Top(u_w)$ to be above $x$, for then $\Top(u_w) \leq \Top(u_{w-1})$. By the lexicographic order of pieces, $u_w$ would be placed in an earlier row than $u_{w-1}$.}
    \label{fig:nonascending}
\end{figure}
Now consider the new row of $P'$ 
inserted between rows $w-1,w$ of $P$.
Let $u_{w-1},u_w$ be the corresponding pieces in $A_c$ 
and $x < y$ be the 
corresponding original rows of $A$.  
It is impossible for $\Top(u_w) = x$ 
because they share the same color.
If $\Top(u_w) < x$ then $\Top(u_{w-1}) \geq \Top(u_w)$. See \Cref{fig:nonascending}(B).
By the lexicographic ordering on the rows of $A_c$, 
$u_w$ would precede $u_{w-1}$, because 
$(\Top(u_w),-\Bottom(u_w)) < (\Top(u_{w-1}),-\Bottom(u_{w-1}))$.
The last case $x < \Top(u_w)$ is also impossible, as
the last 1 in row $\Top(u_w)$ of $A_0$ 
completes the pattern $P'$.
See \Cref{fig:nonascending}(A).
It follows that each $A_c$ is $P$-free.
\[
|A_0| \leq n + \sum_{c=1}^3 |A_c| \leq n + \ex(2n,3m,P).
\]
Taking the $w=1$ case into account, 
we have
$\ex(n,m,P') \leq \ex(2n,3m,P) + n + m$.
\end{proof}

We now turn to \Cref{thm:two-sided-extension}.

\twoextension*

\[
P' = \left(\begin{array}{cccccc} \cline{2-5}
 & \multicolumn{1}{|c}{} &&& \multicolumn{1}{c|}{} &\\
 & \multicolumn{1}{|c}{} && \zero{\hcm[-.32]\rb{2}{$P$}} & \multicolumn{1}{c|}{}\\\cline{2-5}
\bu &&&&& \bu\\\cline{2-5}
 & \multicolumn{1}{|c}{\bu} &&& \multicolumn{1}{c|}{\bu} &\\
  & \multicolumn{1}{|c}{} &&& \multicolumn{1}{c|}{} &\\\cline{2-5}
\end{array}\right).
\]

\begin{proof}
Assume for the time being
that $w>1$.
Let $A_0$ be an $n\times m$, 
$P'$-free matrix.  
Remove the first and last 1 from each row of $A_0$, yielding $A$.
Let $\first(r)$ and $\last(r)$ be the columns
of the first and last 1s 
in row $r$ of $A$.
Associate with row $r$ the 
\emph{bar} $[\first(r) - 1/3, \last(r)+1/3]$.
The bar visibility graph is planar and therefore 4-colorable~\cite{AppelHaken1977,AppelHakenKoch1977,InoueKMMTT26,InoueKMMTT26b}.  Fix any 4-coloring of the rows.
(Extending the bars by 1/3 in both directions
does not change the structure of the visibility graph, but it does ensure that no bar ends at the same position another begins.)

Partition each bar into maximal pieces, 
where each piece has a single upper neighbor.
Obtain $A'$ from $A$ by replacing each 
row $r$ with an interval of rows, one for each piece of $r$,
arranged in ascending or descending order depending on whether $P$ is nonascending or nondescending, respectively.
All pieces in $A'$ originally in row $r$ of $A$ 
inherit the color of $r$.
Define $A_c$ to be the $n_c\times m$ submatrix
of $A'$ on rows colored $c\in \{1,2,3,4\}$.
The first piece in every row of $A$ is charged to the row,
and every subsequent piece is charged to the beginning or ending of some other bar that triggered its creation.  It follows that  $n_1+n_2+n_3+n_4 \leq 3n$.

Suppose $P$ is nondescending and 
$A_c$ contains a copy of $P\in\{0,1\}^{k\times l}$, 
with $P(w,1)=P(w,l)=1$.
We placed the pieces in $A_c$ coming from a 
row $r$ in descending order, which implies that
the embedding of $P$ in $A_c$ maps 
to \emph{distinct} rows in $A$.
Let $x$ and $y$ be the rows in $A$ 
containing pieces $u_{w-1},u_w$
associated with rows $w-1$ and $w$ of $P$.
Let $A(y,j) = A(y,j') = 1$, $j \le j'$, 
be the 1s identified with $P(w,1)=P(w,l)=1$.
Both rows $x,y$ are colored $c$ and therefore nonadjacent.
It follows that $y$'s upper neighbor $\Top(u_w)$
is constant over the interval $[j,j']$ since
it is contained in the single piece $u_w$.
Adding the first and last 1s in row $\Top(u_w)$ of $A_0$
extends the occurrence of $P$ in $A$ to one 
of $P'$ in $A_0$.
\[
P' = 
\begin{array}{r}
\ \\
x\\
\Top(u_w)\\
y\\
\ \\
\end{array}
\left(\begin{array}{cccccc} \cline{2-5}
 & \multicolumn{1}{|c}{} &&& \multicolumn{1}{c|}{} &\\
 & \multicolumn{1}{|c}{} && \zero{\hcm[-.32]\rb{2}{$P$}} & \multicolumn{1}{c|}{}\\\cline{2-5}
\bu &&&&& \bu\\\cline{2-5}
 & \multicolumn{1}{|c}{\bu} &&& \multicolumn{1}{c|}{\bu} &\\
  & \multicolumn{1}{|c}{j} &&& \multicolumn{1}{c|}{j'} &\\\cline{2-5}
\end{array}\right).
\]
The number of 1s in $A_0$ is then
\[
|A_0| \leq 2n + \sum_{c} |A_c| \leq 2n + \sum_c \ex(n_c,m,P)
\leq 2n + \ex(3n,4m,P).
\]
This argument relies on the existence of row $x$, so it must be that $w>1$.  When $w=1$ 
we can pad $P$ with a weight-1 
row 0 with $P(0,1)=1$,
which fixes this problem
but increases the extremal function by at most $m$.
This yields a uniform upper bound 
\[
\ex(n,m,P') \leq \ex(3n,4m,P) + 2n + 4m,
\]
which holds for any $w$.
\end{proof}

\begin{remark}
    \Cref{thm:two-sided-extension} provides another proof that all monotone patterns are linear.  Define $M_k$ to be a monotone pattern on $k$ rows.
    Keszegh~\cite{Keszegh09} proved that $\ex(n,m,M_k)\leq C^k(n+m)$ for some $C>1000$.  Pettie's~\cite{Pettie-SoCG11} grafting lemma
    implies that it is at most $O(5^k(n+m))$.
    \Cref{thm:two-sided-extension} improves this
    to $O(3^k n + 4^k m)$.
\end{remark}

\section{Minimally Nonlinear Patterns}\label{sect:minimally-nonlinear-patterns}

We say a pattern $P$ 
is minimally nonlinear if $\ex(n,P)=\omega(n)$ but
for every $P_0\neq P$ with $P_0\prec P$, $\ex(n,P_0)=O(n)$.\footnote{One could also define minimality to include not just $\prec$ containment but trivial transformations.  
E.g., $Q_3$ would not be minimally nonlinear under this definition, but $Q_3'$ would be.}
In \Cref{sec:hk-minimal} we show that Keszegh's~\cite{Keszegh09} 
family of patterns is minimally nonlinear,
and in \Cref{sec:counting-minimal} we generalize his construction to prove that there are $2^{\Omega(r)}$ minimally nonlinear patterns with weight $r$, improving~\cite{CrowdMath2018,GenesonTsai20}.
\Cref{sect:nlogn-minimally-nonlinear} proves that 
the extremal function of 
any minimally nonlinear pattern is $O(n\log n)$,
except for $C_4$.

\subsection{An Explicit Infinite Family of Minimally Nonlinear Patterns}
\label{sec:hk-minimal}

For $k\geq 0$, let $H_k$ be Keszegh's $(3k+4)\times(3k+4)$ matrix
whose 1s are at positions
\begin{align*}
 (1,2),\ (1,3),\ (3k+2,3k+4),\ (3k+3,3k+4),&\\
 (3j+4,3j+1),
 & \qquad \text{for $0\leq j\leq k$},\\
 b_i=(3i-1,3i+3),\quad c_i=(3i,3i+2), 
 &\qquad \text{for $1\leq i\leq k$}.
\end{align*}
Thus $H_k$ has $3k+5$ ones. We prove the following conjecture of
Keszegh~\cite{Keszegh09}.

\hkminimal*

\Cref{thm:hk-minimal} will follow easily 
from \Cref{thm:hk-minusone-linear}, 
which establishes the linearity of 
another infinite class of matrices.

\begin{theorem}\label{thm:hk-minusone-linear}
    The patterns
    $\bar{P}$ and $P$ are divided into quadrants, and defined below.
\begin{align*}
\bar{P} &= \left(\begin{array}{cc}\bar{A}& Z\\
    0 & \bar{B}
    \end{array}\right),
&&&
P &= \left(\begin{array}{cc}A& Z\\
    0 & B
    \end{array}\right).
\end{align*}
    $\bar{P}$ is a permutation pattern
    whose southwest quadrant is zero
    and whose northeast 
    quadrant $Z$ has a single 1
    not in the top row nor last column.  In $P$, $A$ is obtained from $\bar{A}$ by duplicating the column containing the first-row 1,
    while $B$ is obtained from $\bar{B}$ by duplicating the row containing the last-column 1.
    Then $\ex(n,P) = O(n)$.
\end{theorem}

\begin{proof}
Let $e_{\mathrm{N}}$ and $e_{\mathrm{E}}$ 
(\underline{n}orth and \underline{e}ast)
be the first-row 1 and last-column 1 
in $\bar{P}$.  Let $x$ be the 1 in the column preceding $e_{\mathrm{N}}$ and let $y$ be the 1 in the row below $e_{\mathrm{E}}$.
(If either of $x$ or $y$ does not exist, then $P$ is linear~\cite{Geneson09}.)
Let $z$ be the unique 1 in $Z$.
Let $\bar{A} = 
\left(\begin{array}{c}
\bar{A}_{\mathrm{N}}\\
\bar{A}_{\mathrm{S}}
\end{array}\right)$
and $\bar{B} = \left(\begin{array}{cc}
\bar{B}_{\mathrm{W}} &
\bar{B}_{\mathrm{E}}
\end{array}\right)$
be the partitions of 
$\bar{A},\bar{B}$ 
defined by splitting along 
the row and column 
of $z$, respectively.
Let
$\bar{A} = \left(\begin{array}{cc}
\bar{A}_{\mathrm{W}} &
\bar{A}_{\mathrm{E}}
\end{array}\right)$
be the partition defined by 
cutting between $x$ and $e_{\mathrm{N}}$
and let
$\bar{B} = 
\left(\begin{array}{c}
\bar{B}_{\mathrm{N}}\\
\bar{B}_{\mathrm{S}}
\end{array}\right)$
be the partition defined by cutting
between $y$ and $e_{\mathrm{E}}$.
$A_{\mathrm{N}},A_{\mathrm{S}},\bar{B}_{\mathrm{NW}},\bar{B}_{\mathrm{NE}}$ etc. are defined in the natural way.
A depiction of $P$ is given below, with $x,y,z$ indicated, and the pairs of 1s corresponding to $e_{\mathrm{N}}$ and 
$e_{\mathrm{E}}$ indicated. 
Note that $x$ can be in either $A_{\mathrm{NW}}$ or $A_{\mathrm{SW}}$,
and $y$ can be in either
$B_{\mathrm{SW}}$ or $B_{\mathrm{SE}}$.

\begin{align*}
    P &= \left(
    \begin{array}{cccccccccccccr}
\cline{1-7}
\multicolumn{1}{|c}{} &&\multicolumn{1}{c|}{} & \bu & \bu & & \multicolumn{1}{c|}{} & &&&&&&\\
\multicolumn{1}{|c}{} & \zero{\;A_{\mathrm{NW}}}& \multicolumn{1}{c|}{} & & \zero{\hcm[-.45]\rb{2}{$e_{\mathrm{N}}$}}& & \multicolumn{1}{c|}{} &&&&&&&\\
\multicolumn{1}{|c}{} && \multicolumn{1}{c|}{\zero{\rb{-3}{\hcm[-.1]$x$}}\bu} & & \zero{A_{\mathrm{NE}}} & & \multicolumn{1}{c|}{} &&&&&&&\\
\multicolumn{1}{|c}{} & & \multicolumn{1}{c|}{}
& & & & \multicolumn{1}{c|}{} &&&&&&&\\\cline{1-7}
&&&&&&&&&&\bu\zero{\rb{3}{$z$}}&&&\\\cline{1-7}
\multicolumn{1}{|c}{} && \multicolumn{1}{c|}{} & & & & \multicolumn{1}{c|}{} &&&&&&&\\
\multicolumn{1}{|c}{} &\zero{\;A_{\mathrm{SW}}}& \multicolumn{1}{c|}{} & &\zero{A_{\mathrm{SE}}} & & \multicolumn{1}{c|}{} &&&&&&&\\
\multicolumn{1}{|c}{} && \multicolumn{1}{c|}{}
& & & & \multicolumn{1}{c|}{} &&&&&&&\\\cline{1-10}\cline{12-14}
&&&&&&&\multicolumn{1}{|c}{} &\hcm[.5]& \multicolumn{1}{c|}{} && \multicolumn{1}{|c}{\hcm[.2]} &\istrut{5}\zero{B_{\mathrm{NE}}} & \multicolumn{1}{c|}{\hcm[.2]}\\
&&&&&&&\multicolumn{1}{|c}{} &\zero{\;B_{\mathrm{NW}}} & \multicolumn{1}{c|}{} && \multicolumn{1}{|c}{} && \multicolumn{1}{c|}{\bu}\\
&&&&&&&\multicolumn{1}{|c}{} && \multicolumn{1}{c|}{} && \multicolumn{1}{|c}{} && \multicolumn{1}{c|}{\zero{\hcm[-.5]\rb{3}{$e_{\mathrm{E}}$}}\bu}\\\cline{8-10}\cline{12-14}
&&&&&&&\multicolumn{1}{|c}{} &\zero{\hcm[.85]y}& \multicolumn{1}{c|}{\bu} && \multicolumn{1}{|c}{} && \multicolumn{1}{c|}{}\\
&&&&&&&\multicolumn{1}{|c}{} &\zero{B_{\mathrm{SW}}}& \multicolumn{1}{c|}{} && \multicolumn{1}{|c}{} &\istrut[2]{0}\zero{B_{\mathrm{SE}}}& \multicolumn{1}{c|}{}\\\cline{8-10}\cline{12-14}
    \end{array}
    \right)
\end{align*}

Define $L$ and $R$ as follows,
where $A',B'$ are copies of $A,B$,
respectively.
\begin{align*}
    L &= \left(\;\begin{array}{ccccc}
\cline{1-1}
    \multicolumn{1}{|c|}{\istrut[3]{5}\quad\quad A'_{\mathrm{N}}\quad\quad} & & & & \\\cline{1-2}
    & \multicolumn{1}{|c|}{\istrut[3]{5}\quad\quad A_{\mathrm{N}}\quad\quad} &&&\\\cline{2-2}
    &&&\bu&\\\cline{2-2}
    & \multicolumn{1}{|c|}{\istrut[3]{5}A_{\mathrm{S}}} &&&\\\cline{2-3}\cline{5-5}
    && \multicolumn{1}{|c|}{\istrut[5]{7}\;\bar{B}_{\mathrm{NW}}\;} && \multicolumn{1}{|c|}{\;\bar{B}_{\mathrm{NE}}\;}\\
    &&\multicolumn{1}{|c|}{} && \multicolumn{1}{|c|}{\zero{\rb{-1}{$\quad\quad\;\;\circ$}}}\\
    \cline{1-1}\cline{3-3}\cline{5-5}
    \multicolumn{1}{|c|}{\istrut[3]{5}A'_{\mathrm{S}}} &&&&\\\cline{1-1}\cline{3-3}\cline{5-5}
    && \multicolumn{1}{|c|}{\istrut[6]{9}\bar{B}_{\mathrm{SW}}} && \multicolumn{1}{|c|}{\bar{B}_{\mathrm{SE}}}\\\cline{3-3}\cline{5-5}
    \end{array}\;\right),
&&&\\
&&&\\
R &= \left(\;\begin{array}{ccccccc}
\cline{1-1}\cline{3-3}
\multicolumn{1}{|c|}{\istrut[6]{9}\;\bar{A}_{\mathrm{NW}}\;}& & 
\multicolumn{1}{|c|}{\zero{\hcm[-0.3]\rb{6.5}{$\;\;\circ$}}\;\bar{A}_{\mathrm{NE}}\;} &&&&\\\cline{1-1}\cline{3-3}
&&&&\bu&&\\\cline{1-1}\cline{3-3}
\multicolumn{1}{|c|}{\istrut[6]{9}\; \bar{A}_{\mathrm{SW}}\;} & & 
\multicolumn{1}{|c|}{\;\bar{A}_{\mathrm{SE}}\;} &&&&\\\cline{1-1}\cline{3-4}\cline{6-6}
&&& \multicolumn{1}{|c|}{\istrut[3]{5}\quad B_{\mathrm{W}}\quad} & &
\multicolumn{1}{|c|}{\quad B_{\mathrm{E}}\quad}\\\cline{2-2}\cline{4-4}\cline{6-7}
& \multicolumn{1}{|c|}{\istrut[3]{5}\quad B'_{\mathrm{W}}\quad} &&&&& 
\multicolumn{1}{|c|}{\quad B'_{\mathrm{E}}\quad}\\\cline{2-2}\cline{7-7}
\end{array}\;\right).
\end{align*}
The 1 in the southeast corner of $\bar{B}_{\mathrm{NE}}$ corresponds 
to $e_{\mathrm{E}}$
and the 1 in the northwest corner
of $\bar{A}_{\mathrm{NE}}$ corresponds to $e_{\mathrm{N}}$,
both indicated by a $\circ$.

Observe that $L$ is contained in a double permutation pattern.  
Except for the repeated columns in $A_{\mathrm{N}}'$ and $A_{\mathrm{N}}$,
all 1s appear in distinct rows and columns.
Similarly $R$ is contained in the rotation 
of a double permutation.
According to 
Geneson's bound~\cite{Geneson09},
$\ex(n,L), \ex(n,R) = O(n)$.

Let $A$ be a $P$-free matrix.  
Assign the label \textsf{L} to any 1 in $A$
that plays the role of $e_{\mathrm{E}}$ in an occurrence of $L$, and the label 
\textsf{R} to any 1 in $A$
that plays the role of $e_{\mathrm{N}}$ in $R$.
The part of $A$ not labeled \textsf{L} is $L$-free 
and the part not labeled \textsf{R} is $R$-free.
Therefore, if every 1 has at most one label
then $|A| \leq \ex(n,L) + \ex(n,R) = O(n)$.
We claim that this is true, and for the purpose of obtaining a contradiction suppose there exists a 
1 in $A$ labeled both \textsf{L} and \textsf{R},
which is named $e$ in \Cref{fig:large-LR-diagram}.

\begin{figure}[h]
\centering
\input{large-diagram}
\caption{\label{fig:large-LR-diagram}The 1 labeled $e$ is identified with occurrences of 
$e_{\mathrm{E}}$ in $L$ 
and
$e_{\mathrm{N}}$ in $R$, respectively.}
\end{figure}
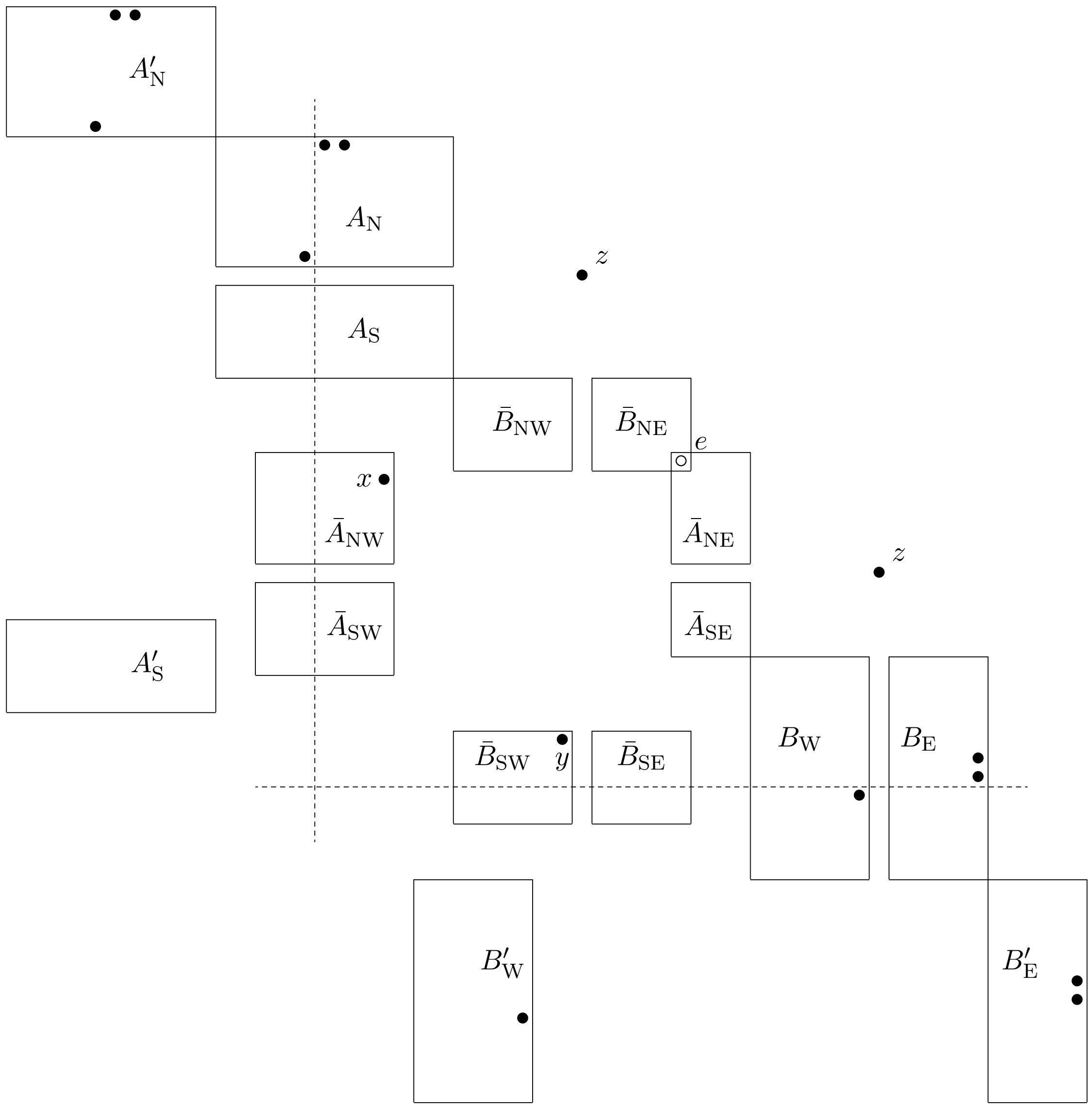

Let ``$x$'' and ``$y$'' refer to their copies in $\bar{A}$ (in $R$) and $\bar{B}$ (in $L$), 
respectively.
If we draw a horizontal line just below the 
repeated row 
in $B_{\mathrm{E}}$ then $y$
must be above this line, for otherwise 
$P\prec A$, as witnessed by 
$A_N,A_S,\bar{B}_{NW}, \bar{B}_{NE},\bar{B}_{SW}, \bar{B}_{SE}$,
the $z$ from $L$,
and the last column of $B_E$.
Similarly, if we draw a vertical line just to the 
left of the repeated column in $A_N$,
$x$ must be to the right of this line,
for otherwise $P\prec A$ as witnessed
by $\bar{A}_{NW},\bar{A}_{NE},\bar{A}_{SW},\bar{A}_{SE},B_W,B_E$, the $z$ from $R$,
and the top row of $A_N$.
On the other hand, if $x$ is to the 
right of its line then so is $B_{W}'$,
and if $y$ is above its line then so is $A_S'$.
In this case we also have 
$P\prec A$, as witnessed by
$A_N',A_S',B_W',B_E',$ and $e$, which plays 
the role of $z$.

All possibilities lead to contradiction, 
so we conclude that 
$\ex(n,P) \leq \ex(n,L) + \ex(n,R) = O(n)$.
\end{proof}

\Cref{thm:hk-minusone-linear} suffices to 
prove \Cref{thm:hk-minimal}.

\begin{proof}[Proof of Theorem~\ref{thm:hk-minimal}]
Divide $H_k$ into quadrants by cutting between columns $3i$ and $3i+1$ and between 
rows $3i+1$ and $3i+2$.  
By definition of $H_k$,
the southwest quadrant is zero, 
and the northeast quadrant contains only
$b_i$ and $c_i$.  
The top row contains 1s in the second and third columns, and the last column contains 1s in the second-to-last and third-to-last rows.
Except for these duplicates, every 1 is 
isolated in its row and column.
\begin{align*}
    \left(\;\begin{array}{cccccccc}\cline{1-4}
    \multicolumn{1}{|c}{} & \bu & \bu & \multicolumn{1}{c|}{} &&&&\\
    \multicolumn{1}{|c}{} & & & \multicolumn{1}{c|}{} &&&&\\
    \multicolumn{1}{|c}{} &  & & \multicolumn{1}{c|}{} &&\zero{\rb{3}{$\;\:b_i$}}\bu&&\\
    \multicolumn{1}{|c}{} & & & \multicolumn{1}{c|}{} &\zero{\rb{3}{$\;\:c_i$}}\bu&&&\\\cline{1-8}
    &&&& \multicolumn{1}{|c}{} &&& \multicolumn{1}{c|}{} \\
    &&&& \multicolumn{1}{|c}{} &&& \multicolumn{1}{c|}{\bu} \\
    &&&& \multicolumn{1}{|c}{} &&& \multicolumn{1}{c|}{\bu} \\
    &&&& \multicolumn{1}{|c}{} &&& \multicolumn{1}{c|}{} \\\cline{5-8}
    \end{array}\;\right)
\end{align*}
If we remove either $b_i$ or $c_i$ the resulting pattern falls within the class covered
by \Cref{thm:hk-minusone-linear} and is therefore linear.  Removing one of the duplicate 1s in the first row or last column leaves a double permutation matrix, which is linear by Geneson~\cite{Geneson09}. 
Removing any of the other 1s, 
at positions $(3j+4,3j+1)$, $0\leq j\leq k$,
leaves a matrix known to be linear
from \cite{Geneson09} and Keszegh's~\cite{Keszegh09} 
gluing operation.
Since $H_k$ is nonlinear~\cite{Keszegh09,Pettie-FH11}, it is \emph{minimally} nonlinear.
\end{proof}

\subsection{Exponentially many minimally nonlinear matrices}
\label{sec:counting-minimal}

CrowdMath~\cite{CrowdMath2018} and Geneson and Tsai~\cite{GenesonTsai20} proved upper bounds of $r^{O(r)}$ for the number of minimally nonlinear patterns with $r$ rows. We give an exponential lower bound, even when the number of columns and the number of 1s are fixed.

\begin{theorem}\label{thm:exponentially-many-minimal}
For all sufficiently large $r$, 
there are $2^{\Omega(r)}$ minimally nonlinear $r\times r$ matrices with weight $r+1$.
Each matrix $P$ in the construction satisfies $\ex(n,P)=\Theta(n\log n)$, 
where the constants may depend on $P$.
\end{theorem}

We use the following lower bound from 
Pettie~\cite[Theorem~3.1]{Pettie-FH11}, 
fixing the permutation in that construction 
to be the identity.

\begin{theorem}[{Pettie~\cite{Pettie-FH11}}]\label{thm:pettie-increasing-cut}
Let $Q$ be a $0$--$1$ matrix with at least two 1s. 
Suppose that it is impossible to partition $Q$ 
into quadrants, some of which may be empty, 
satisfying the 
following conditions.
\begin{enumerate}
\item The southwest quadrant contains no 1s.
\item The 1s in the northeast quadrant are contained in a sufficiently large identity matrix $I$.
\item At least two quadrants are nonzero.
\end{enumerate}
Then $\ex(n,Q)=\Omega(n\log n)$.
\end{theorem}

\begin{proof}[Proof of \Cref{thm:exponentially-many-minimal}]
We construct a minimally nonlinear pattern as follows.  Pick integers $q$ and $\ell \leq q$,
and a composition 
of $q=m_1+\cdots+m_\ell$ into $\ell$ positive parts.
Let $I_1,\ldots,I_\ell$ be the corresponding partition of $[q]$ into intervals,
with $|I_i|=m_i$.

Begin with the following ordered 
row and column labels.
\[
\begin{array}{rl}
\text{rows:} &
H, \;p_1, \;z_0, \;p_2, \;z_1, \ldots, \;p_q, \;z_{q-1}, \;V_1, \;V_2, \;z_q,\\
\text{columns:} &
z_0, \;H_1, \;H_2, \;z_1, \;p_1, \;z_2, \;p_2,\ldots, \;z_q, \;p_q, \;V.
\end{array}
\]
For each interval $I_j$, insert a row labeled $t_j$ immediately before its first $p$ row and a column labeled $t_j$ immediately after its last $p$ column. The pattern $P = P(m_1,\ldots,m_\ell)$ 
is now formed as follows.
Put 1s at all positions with matching row and column labels, in row $H$ at columns $H_1,H_2$, and in column $V$ at rows $V_1,V_2$. 
All other entries of $P$ are zero. 
We also use $p_i,t_i,z_i$ to refer to the 1s 
at the positions indicated by those 
row and column labels.

The pattern $P$ has $r=2q+\ell+4$ rows and columns and $r+1$ 1s. Its first row has 1s in columns $2,3$, and its last column has 1s in rows $r-2,r-1$; all other rows and columns have a single 1. When $m_1,\ldots,m_q=1$, $P$ is Keszegh's 
$H_q$ from \Cref{sec:hk-minimal}. \Cref{fig:composition-patterns} illustrates 
the construction of $P(1,2)$ and $P(2,1)$.

\begin{figure}[htbp]
\centering
\begin{tikzpicture}[x=.43cm,y=.43cm,baseline=(current bounding box.center)]
  \foreach \i in {0,...,12} {
    \draw[gray!25,thin] ({\i-.5},.5) -- ({\i-.5},-11.5);
    \draw[gray!25,thin] (-.5,{.5-\i}) -- (11.5,{.5-\i});
  }
  \draw[gray!60,thin] (-.5,-11.5) rectangle (11.5,.5);
  \node[anchor=east,font=\footnotesize] at (-.7,0) {$H$};
  \node[anchor=east,font=\footnotesize] at (-.7,-1) {$t_1$};
  \node[anchor=east,font=\footnotesize] at (-.7,-2) {$p_1$};
  \node[anchor=east,font=\footnotesize] at (-.7,-3) {$z_0$};
  \node[anchor=east,font=\footnotesize] at (-.7,-4) {$t_2$};
  \node[anchor=east,font=\footnotesize] at (-.7,-5) {$p_2$};
  \node[anchor=east,font=\footnotesize] at (-.7,-6) {$z_1$};
  \node[anchor=east,font=\footnotesize] at (-.7,-7) {$p_3$};
  \node[anchor=east,font=\footnotesize] at (-.7,-8) {$z_2$};
  \node[anchor=east,font=\footnotesize] at (-.7,-9) {$V_1$};
  \node[anchor=east,font=\footnotesize] at (-.7,-10) {$V_2$};
  \node[anchor=east,font=\footnotesize] at (-.7,-11) {$z_3$};
  \node[anchor=west,rotate=60,font=\footnotesize] at (0,.7) {$z_0$};
  \node[anchor=west,rotate=60,font=\footnotesize] at (1,.7) {$H_1$};
  \node[anchor=west,rotate=60,font=\footnotesize] at (2,.7) {$H_2$};
  \node[anchor=west,rotate=60,font=\footnotesize] at (3,.7) {$z_1$};
  \node[anchor=west,rotate=60,font=\footnotesize] at (4,.7) {$p_1$};
  \node[anchor=west,rotate=60,font=\footnotesize] at (5,.7) {$t_1$};
  \node[anchor=west,rotate=60,font=\footnotesize] at (6,.7) {$z_2$};
  \node[anchor=west,rotate=60,font=\footnotesize] at (7,.7) {$p_2$};
  \node[anchor=west,rotate=60,font=\footnotesize] at (8,.7) {$z_3$};
  \node[anchor=west,rotate=60,font=\footnotesize] at (9,.7) {$p_3$};
  \node[anchor=west,rotate=60,font=\footnotesize] at (10,.7) {$t_2$};
  \node[anchor=west,rotate=60,font=\footnotesize] at (11,.7) {$V$};
  \fill (1,0) circle (1.65pt);
  \fill (2,0) circle (1.65pt);
  \fill (5,-1) circle (1.65pt);
  \fill (4,-2) circle (1.65pt);
  \fill (0,-3) circle (1.65pt);
  \fill (10,-4) circle (1.65pt);
  \fill (7,-5) circle (1.65pt);
  \fill (3,-6) circle (1.65pt);
  \fill (9,-7) circle (1.65pt);
  \fill (6,-8) circle (1.65pt);
  \fill (11,-9) circle (1.65pt);
  \fill (11,-10) circle (1.65pt);
  \fill (8,-11) circle (1.65pt);
  \draw[black!70,densely dashed,line width=.45pt] (-.5,-6.5) -- (11.5,-6.5);
  \draw[black!70,densely dashed,line width=.45pt] (5.5,.5) -- (5.5,-11.5);
  \node[anchor=north,font=\small] at (5.5,-12.1) {$P(1,2)$};
\end{tikzpicture}
\hspace{1.1cm}
\begin{tikzpicture}[x=.43cm,y=.43cm,baseline=(current bounding box.center)]
  \foreach \i in {0,...,12} {
    \draw[gray!25,thin] ({\i-.5},.5) -- ({\i-.5},-11.5);
    \draw[gray!25,thin] (-.5,{.5-\i}) -- (11.5,{.5-\i});
  }
  \draw[gray!60,thin] (-.5,-11.5) rectangle (11.5,.5);
  \node[anchor=east,font=\footnotesize] at (-.7,0) {$H$};
  \node[anchor=east,font=\footnotesize] at (-.7,-1) {$t_1$};
  \node[anchor=east,font=\footnotesize] at (-.7,-2) {$p_1$};
  \node[anchor=east,font=\footnotesize] at (-.7,-3) {$z_0$};
  \node[anchor=east,font=\footnotesize] at (-.7,-4) {$p_2$};
  \node[anchor=east,font=\footnotesize] at (-.7,-5) {$z_1$};
  \node[anchor=east,font=\footnotesize] at (-.7,-6) {$t_2$};
  \node[anchor=east,font=\footnotesize] at (-.7,-7) {$p_3$};
  \node[anchor=east,font=\footnotesize] at (-.7,-8) {$z_2$};
  \node[anchor=east,font=\footnotesize] at (-.7,-9) {$V_1$};
  \node[anchor=east,font=\footnotesize] at (-.7,-10) {$V_2$};
  \node[anchor=east,font=\footnotesize] at (-.7,-11) {$z_3$};
  \node[anchor=west,rotate=60,font=\footnotesize] at (0,.7) {$z_0$};
  \node[anchor=west,rotate=60,font=\footnotesize] at (1,.7) {$H_1$};
  \node[anchor=west,rotate=60,font=\footnotesize] at (2,.7) {$H_2$};
  \node[anchor=west,rotate=60,font=\footnotesize] at (3,.7) {$z_1$};
  \node[anchor=west,rotate=60,font=\footnotesize] at (4,.7) {$p_1$};
  \node[anchor=west,rotate=60,font=\footnotesize] at (5,.7) {$z_2$};
  \node[anchor=west,rotate=60,font=\footnotesize] at (6,.7) {$p_2$};
  \node[anchor=west,rotate=60,font=\footnotesize] at (7,.7) {$t_1$};
  \node[anchor=west,rotate=60,font=\footnotesize] at (8,.7) {$z_3$};
  \node[anchor=west,rotate=60,font=\footnotesize] at (9,.7) {$p_3$};
  \node[anchor=west,rotate=60,font=\footnotesize] at (10,.7) {$t_2$};
  \node[anchor=west,rotate=60,font=\footnotesize] at (11,.7) {$V$};
  \fill (1,0) circle (1.65pt);
  \fill (2,0) circle (1.65pt);
  \fill (7,-1) circle (1.65pt);
  \fill (4,-2) circle (1.65pt);
  \fill (0,-3) circle (1.65pt);
  \fill (6,-4) circle (1.65pt);
  \fill (3,-5) circle (1.65pt);
  \fill (10,-6) circle (1.65pt);
  \fill (9,-7) circle (1.65pt);
  \fill (5,-8) circle (1.65pt);
  \fill (11,-9) circle (1.65pt);
  \fill (11,-10) circle (1.65pt);
  \fill (8,-11) circle (1.65pt);
  \draw[black!70,densely dashed,line width=.45pt] (-.5,-3.5) -- (11.5,-3.5);
  \draw[black!70,densely dashed,line width=.45pt] (4.5,.5) -- (4.5,-11.5);
  \node[anchor=north,font=\small] at (5.5,-12.1) {$P(2,1)$};
\end{tikzpicture}
\caption{Two $12\times12$ matrices with thirteen 1s. The dashed cuts leave $p_2,t_2$ northeast and no 1s southwest on the left, and $t_1$ northeast and $z_1$ southwest on the right.}
\label{fig:composition-patterns}
\end{figure}
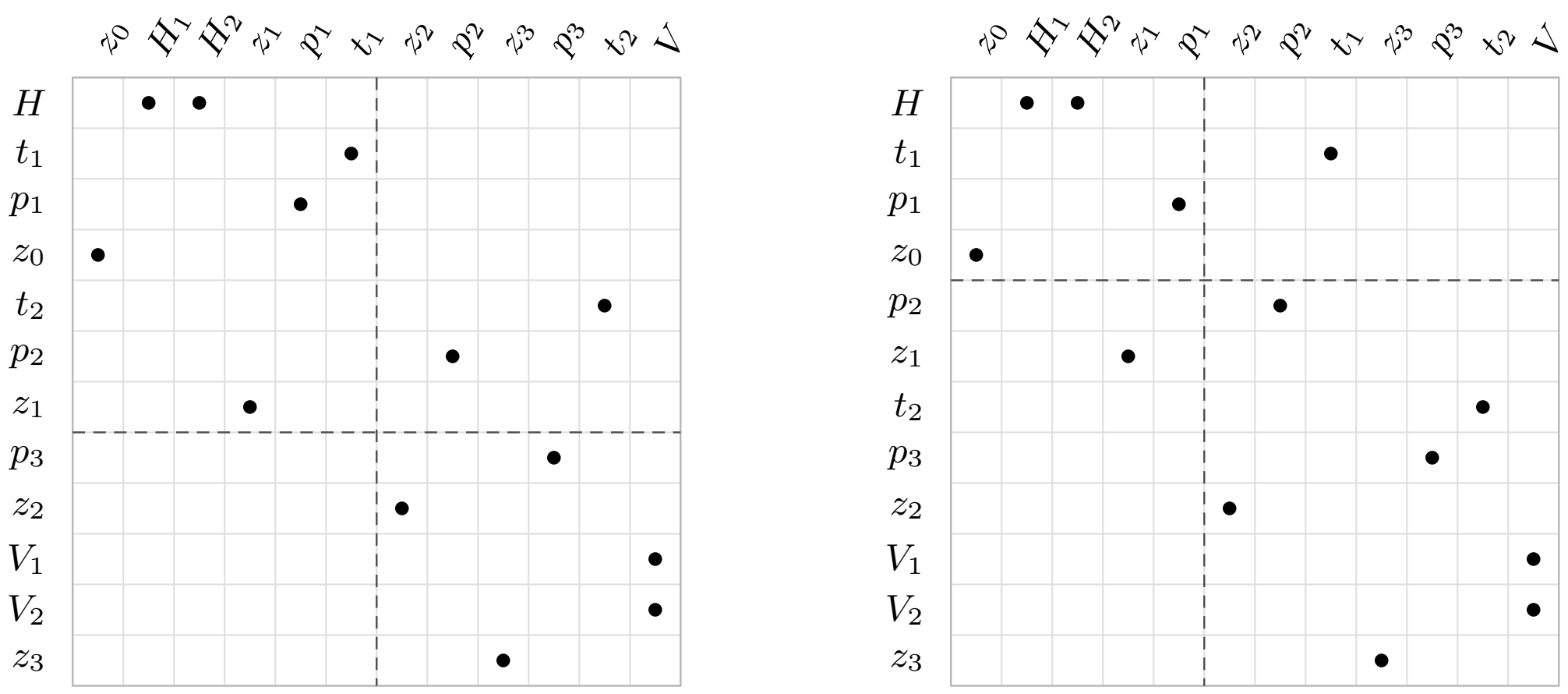

For $i\in I_j$, consider a horizontal cut after 
row $z_{i-1}$ and either of the following 
vertical cuts:
\begin{center}
\begin{tabular}{l|c|c}
\multicolumn{1}{l}{\textsf{The Vertical Cut}} & \multicolumn{1}{c}{\textsf{NE Quadrant 1s}} & 
\multicolumn{1}{c}{\textsf{SW Quadrant 1s}}\\\hline
\istrut[2]{4}Immediately before column $z_i$ 
& $p_i,t_j$ & $\emptyset$\\\hline
\istrut[2]{4}Immediately after column $p_i$ & $t_j$ & $z_i$.\\\hline
\end{tabular}
\end{center}
The 1s above the horizontal cut are exactly the two 1s in row $H$, the entries $p_a$ with $1\leq a\leq i$, the entries $z_a$ with $0\leq a<i$, and the entries $t_b$ with $1\leq b\leq j$. Among these, only $p_i,t_j$ lie to the right of the first cut, since the earlier blocks end before $i$ and the current block ends at or after $i$. Every later 1 lies to its right as well, so in this case 
the SW quadrant is zero.

Moving the cut past the consecutive columns $z_i,p_i$ gives the second line of the table. Note that $p_i,t_j$ are in decreasing order: $t_j$ has the earlier row and later column.

\paragraph{Nonlinearity.}
By \Cref{thm:pettie-increasing-cut}, it suffices to rule out cuts satisfying the  stated conditions. Suppose such cuts exist. Since the northeast quadrant is nonempty, row $H$ lies above the horizontal cut. 
We consider the first $z$ entry below the horizontal cut, if there is one.
\begin{itemize}
\item If it is $z_0$, the zero southwest quadrant forces the vertical cut before the first column, putting both $H$ 1s in the northeast quadrant.
\item If it is $z_i$ with $i>0$, the horizontal cut is after $z_{i-1}$ and the vertical cut is before $z_i$. Thus the northeast quadrant contains the decreasing pair $p_i,t_j$, where $i\in I_j$.
\item If no $z$ entry is below the cut, all rows are above it, and a nonempty northeast quadrant contains both 1s in the last column.
\end{itemize}
These cases are exhaustive, 
and certify that there is no partition of $P$ into quadrants that satisfies the criteria 
of \Cref{thm:pettie-increasing-cut}.
Thus $\ex(n,P)=\Omega(n\log n)$.

\paragraph{Minimality.}
We show that deleting any 1-entry from $P$
leaves a linear pattern. 
\begin{enumerate}
\item \emph{Deleting a 1 in $H$ or $V$.} 
Deleting a 1 from $H$ or $V$ leaves a pattern that is contained in a double permutation or rotation thereof, all of which are linear~\cite{Geneson09}.

\item \emph{Deleting $z_0$ or $z_q$.} 
After deleting $z_0$, deleting the 1 in column $H_1$ reduces the extremal function by at most $n$.  Similarly, after deleting $z_q$, deleting the 1 in row $V_2$ reduces the extremal function by at most $n$.  In both cases, what remains is contained in a double permutation pattern, which is linear~\cite{Geneson09}.

\item \emph{All remaining 1s.} 
Fix $j$ and $i\in I_j$. 
After deleting $p_i$ or $t_j$, 
use the first cut in the table.
After deleting $z_i$, use the second cut. 
Each cut leaves exactly one isolated 1 in the 
northeast quadrant and an all-zero southwest quadrant. 
The resulting pattern is a permutation,
except for the double 1s in the first row of the northwest quadrant, and the double 1s in the last column of the southeast quadrant.
\Cref{thm:hk-minusone-linear} 
states that such patterns are linear.
\end{enumerate}

Deleting any 1 leaves a linear pattern so 
$P$ is minimally nonlinear.  Furthermore, 
since $P\neq C_4$, \Cref{thm:minimally-nonlinear} implies an upper bound of $\ex(n,P)=O(n\log n)$.

\paragraph{Counting.}
Recall that $P$ is an $r\times r$ pattern with weight $r+1$,
where $r = 2q+\ell+4$.
Define $a_s$ to be the number of 
compositions $(m_1,\ldots,m_\ell)$
for which each $m_i\geq 1$, 
$\sum_i m_i = q$ and $2q+\ell = s$.  
One may confirm that 
$a_5 = a_6 = a_7 = 1$.
For $s\geq 8$, we can take 
any eligible vector with $2q+\ell=s$ 
and either decrement $m_\ell$ if it is greater than 1, or remove it if it is 
equal to 1.  This mapping gives a bijection between compositions counted by $a_s$ and those counted by $a_{s-3}$ and $a_{s-2}$, therefore
\[
a_s = a_{s-2} + a_{s-3}.
\]
It follows that $a_s = 2^{\Omega(s)}$
and that there are $2^{\Omega(r)}$ minimally nonlinear $r\times r$ patterns with weight $r+1$.
\end{proof}

\subsection{Extremal Functions of Minimally Nonlinear Matrices}\label{sect:nlogn-minimally-nonlinear}

In this section we prove \Cref{thm:minimally-nonlinear}.

\minimallynonlinear*

We first prove a technical lemma.  Recall the 
patterns $C_4,R_1,R_3,R_4$.
\begin{align*}
C_4 &= \patbits{11/11},
&R_1 &= \patbits{110/101},
&R_3 &= \patbits{110/100/001},
&R_4 &= \patbits{1000/0001/0001/0110}.
\end{align*}

\begin{lemma}\label{lem:boundary-singleton}
Let $P$ be a pattern that avoids all rotations and reflections of
$R_1,R_3,R_4,C_4$, 
and has no all-zero rows or columns.
Then among the first and last rows of $P$, and the first and last columns of $P$, 
at least one has weight 1.
\end{lemma}

\begin{proof}
Let $P\in \{0,1\}^{a\times b}$ be the pattern and suppose, 
for the purpose of obtaining a contradiction, that all exterior rows/columns of $P$ have weight at least 2.

\paragraph{Case: two adjacent corners.} Suppose 
$P(1,1) = P(1,b)=1$. By assumption there exist $i,i'>1$ 
such that $P(i,1)=P(i',b)=1$,
which forms $C_4$ if $i=i'$
and a rotation/reflection of $R_1$ if $i\neq i'$.

\paragraph{Case: two opposite corners.}
Suppose $P(1,1)=P(a,b)=1$.
The first row and first column each have another 1, which gives an occurrence of $R_3$.

\paragraph{Case: one corner.}
Suppose $P(1,1)=1$, but all other corners are zero. The last row and column each have two 1s, which gives an occurrence of $R_4$.

\paragraph{Case: no corners.}
When all corners are zero,
let $j,j' \in (1,b)$ 
be such that 
$P(1,j)=P(1,j')=1$, 
and let $i_0 > i_0' > 1$ and
$i_1 > i_1' > 1$ be such that 
$P(i_0',1)=P(i_0,1)=P(i_1',b)=P(i_1,b)=1$.
If $i_0\neq i_1$ we have a rotation/reflection 
of $R_4$ in $P$, and if $i_0=i_1$
we have a rotation/reflection of $R_1$ in $P$.
\[
\begin{array}{r}
\ \\
i_1'\\
\ \\
i_1\\
i_0\\
\end{array}
\patbits{01010/00000/00001/00000/00001/10000/00000},
\qquad
\begin{array}{r}
\ \\
i_0'\\
i_1'\\
\ \\
i_0\\
\end{array}
\patbits{01010/00000/10000/00001/00000/10001/00000}.
\]
\end{proof}

\medskip 

\begin{proof}[Proof of \Cref{thm:minimally-nonlinear}]
If $P$ contains a rotation/reflection of 
$R_1,R_3,R_4,$ or $C_4$ then by minimality
$P$ must be equal to one of them.
\Cref{lem:boundary-singleton} implies
that one of $P$'s exterior rows/columns has 
weight 1.  Delete that row/column to 
obtain a proper subpattern $P'\prec P$. 
Minimality of $P$ implies 
$\ex(n,P')=O(n)$, while the Pach--Tardos transformation implies
$\ex(n,P) = O(\ex(n,P')\log n)=O(n\log n)$.
\end{proof}

\section{Patterns with Extremal Function $\Theta(n\log n/\log\log n)$}\label{sect:logn-over-loglogn}

In this section we give \emph{general} methods for proving $\Omega(n\log n/\log\log n)$ lower bounds by generalizing Tardos's construction~\cite{Tardos05}, and $O(n\log n/\log\log n)$ upper bounds using a simple divide and conquer strategy.
We apply our methods to the weight-6 forbidden patterns $S_1,\ldots,S_{12}$, 
obtaining sharp $\Theta(n\log n/\log\log n)$ bounds on 
$S_1,\ldots,S_{10}$.\footnote{$S_1,\ldots,S_{12}$ are exhaustive of weight-6 patterns susceptible to our lower bounds.  There are many other weight-6 patterns that can be upper bounded by $O(n\log n/\log\log n)$ via \Cref{thm:logn-over-loglogn-decomposition}.}
Several of the upper bounds depend on
\Cref{thm:Qhat-alpha}, 
proved in Appendix~\ref{sect:new-quasilinear-bounds},
which states that $\ex(n,m,\hat{Q}_3) = O(n\alpha(n,m)+m)$.

\begin{align*}
\begin{array}{lll}
    S_1 = \patbits{00110/01001/00100/10000}, &
    S_2 = \patbits{00110/01000/00001/00100/10000}, &
    S_3 = \patbits{01010/10100/00101},\\&\\
    S_4 = \patbits{01010/00100/10000/00101}, &
    S_5 = \patbits{01010/00101/00100/10000}, &
    S_6 = \patbits{01010/00100/00001/00100/10000},\\&\\
    S_7 = \patbits{01010/00100/10101}, &
    S_8 = \patbits{01010/00100/00101/10000}, &
    S_9 = \patbits{01010/00100/00100/10001}, \\&\\
    S_{10} = \patbits{01010/00100/00100/00001/10000}, &
    S_{11} = \patbits{0110/1000/0001/0100/1000},&
    S_{12} = \patbits{1010/0001/0010/0100/0001}.
\end{array}
\end{align*}

\subsection{Lower Bounds}\label{subsect:logn-over-loglogn-lb}

\subsubsection{Tardos's Construction}\label{subsubsect:Tardos-construction}

We use Tardos's original 
construction~\cite{Tardos05} and variations on it.
$A_{k,r}$ is a $kr^{k-1} \times r^{k}$ 0--1 matrix, 
whose columns are identified with vectors in 
$\{1,2,\ldots,r\}^k$ and whose rows are 
identified with those vectors in $\{0,1,2,\ldots,r\}^k$ 
having exactly one 0.
\[
A(u,v) = \begin{cases}
    1 & \text{ if $v-u$ has one non-zero entry,}\\
    0 & \text{ otherwise.}
\end{cases}
\]
In other words, $v$ is derived from $u$ by replacing its 
0 with a digit in $\{1,2,\ldots,r\}$.  
$A_{k,r}$ can also be defined recursively as follows,
where $I$ is the $r^{k-1}\times r^{k-1}$ identity matrix.
\[
A_{k,r} \,= \;\:
\begin{array}{rccccc}
&1&2&3&\cdots&r\\\cline{2-6}
0 & \multicolumn{1}{|c}{\hcm[.5]I\hcm[.5]\istrut[2]{5}}
& \multicolumn{1}{|c}{\hcm[.5]I\hcm[.5]}
& \multicolumn{1}{|c}{\hcm[.5]I\hcm[.5]}
& \multicolumn{1}{|c}{\hcm[.5]I\hcm[.5]}
& \multicolumn{1}{|c|}{\hcm[.5]I\hcm[.5]}\\\cline{2-6}
1 & \multicolumn{1}{|c|}{A_{k-1,r}\istrut[2]{5}}
& \multicolumn{4}{c|}{}\\\cline{2-3}
\multicolumn{1}{r|}{2} &\multicolumn{1}{c|}{}
& \multicolumn{1}{c|}{A_{k-1,r}\istrut[2]{5}}
& \multicolumn{3}{c|}{}\\\cline{3-4}
\multicolumn{1}{r|}{3} &\multicolumn{2}{c|}{}
& \multicolumn{1}{c|}{A_{k-1,r}\istrut[2]{5}}
& \multicolumn{2}{c|}{}\\\cline{4-5}
\multicolumn{1}{r|}{\vdots\hcm[.04]} &\multicolumn{3}{c|}{}
& \multicolumn{1}{c|}{\ddots\istrut[2]{5}}
& \multicolumn{1}{c|}{}\\\cline{5-6}
\multicolumn{1}{r|}{r} 
& \multicolumn{4}{c|}{}
& \multicolumn{1}{c|}{A_{k-1,r}\istrut[2]{5}}\\\cline{2-6}
\end{array}
\]

The rows and columns are labeled by the first digit of their vectors.  Consider the row-0 blocks.  When $u = 0\hat{u}$, it is incident to exactly one 1 in each of the row-0 blocks with column labels $1,2,\ldots,r$.  
Because both rows and columns are ordered lexicographically, 
all row-0 blocks are the identity $I$. 
The row-$w$ blocks with $w\neq 0$ as their first digit can only have 1s at the intersection of row-$w$ and column-$w$.  It also follows that each of these blocks is a copy of $A_{k-1,r}$.  

A small variation on this recursive construction is to substitute for $I$ the $90^\circ$ rotation of $I$.
Call it $I^{\text{rot}}$ and 
the resulting matrix $A_{k,r}^{\text{rot}}$.

\begin{theorem}\label{thm:lex-lex-lowerbound}
    Let $P\in\{0,1\}^{a\times b}$ be a pattern with $a b > 1$
    and $r > \max\{a,b\}$.
    Define $B$ to be a 0--1 matrix 
    with rows and columns partitioned into
    intervals labeled $0,1,2,\ldots,r$ and $1,2,\ldots,r$, respectively, whose intersections form \emph{blocks}.  
    Block $(0,t) = I$ and block 
    $(t,t)$ is the all-1 matrix, 
    for $t\in \{1,2,\ldots,r\}$.
    The dimensions of the blocks are at least $a+b$.
    Define $B^{\mathrm{rot}}$ analogously, using $I^{\mathrm{rot}}$ in lieu of $I$.
    If there is no nontrivial embedding of 
    $P$ in $B$, or no nontrivial embedding of $P$ in $B^{\mathrm{rot}}$
    then $\ex(n,P) = \Omega(n\log n/\log\log n)$.
    An embedding is trivial if it puts $P$ entirely in one block.
\end{theorem}

\begin{proof}
    When $k\ge r$, $A = A_{k,k}$ is a square matrix with $k^{k+1}$ 1s. Setting $n=k^k$, we have $|A| = \Theta(n\log n/\log\log n)$.
    If $P\nprec A$ then $\ex(n,P)=\Omega(n\log n/\log\log n)$.  Suppose $P\prec A$, and let $t$ be the
    \emph{first} digit where any of the rows or columns involved in the embedding of $P$ in $A$ differ.  
    Then $P$ has a \emph{nontrivial} embedding in $A_{k-t+1,k}$, 
    and hence a nontrivial embedding in $B$ as well, contradicting the assumption on $P$.
\end{proof}

\begin{corollary}\label{cor:S-one-through-ten}
For $k\in [11]$, $\ex(n,S_k) = \Omega(n\log n/\log\log n)$.
In particular, 
\begin{align*}
S_1,S_2,S_3,S_4,S_7,S_8,S_9,S_{10},S_{11} &\nprec A,\\
S_5,S_6,S_7,S_8,S_9,S_{10} &\nprec A^{\mathrm{rot}}.
\end{align*}
\end{corollary}

\begin{proof}
    We illustrate the argument for $S_1$.
    Partition the nonzero blocks of $B$ into the \emph{row-0 blocks} and the \emph{diagonal blocks}.  $S_1$ cannot be nontrivially embedded in just the diagonal blocks because it is not block diagonal, so at least one row must be embedded in the row-0 blocks.  If \emph{only} the first row of $S_1$ 
    is embedded in the row-0 blocks
    then the two 1s in the top row must appear in different row-0 blocks, and the remaining rows cannot be embedded in the diagonal blocks.  If only the first two rows are embedded in the row-0 blocks, then three 1s in the first two rows of $S_1$ must be in distinct row-0 blocks, leaving the remaining rows impossible to embed in the diagonal blocks.  These cases are illustrated below.
    \[
    \left(\begin{array}{ccc|cc}
    &&\bu&\bu&\\\hline
    &\bu&&&\bu\\
    &&\bu&&\\
    \bu&&&&\end{array}
    \right),
    \qquad
    \left(\begin{array}{cc|c|cc}
    &&\bu&\bu&\\
    &\bu&&&\bu\\\hline
    &&\bu&&\\
    \bu&&&&\end{array}
    \right).
    \]
\end{proof}

\subsubsection{Tardos's Construction, Antilexicographically}\label{subsubsect:antilex-Tardos}

Define $A^{\text{anti}}$ to be $A$, 
but with columns ordered anti-lexicographically,
that is, they appear according to 
the lexicographic order of their reversals.
The definition of a \emph{legal division} captures the properties necessary to embed a pattern in $A^{\text{anti}}$.

\begin{definition}[Span and Legal Division]\label{def:legal-division}
    Let $P\in\{0,1\}^{k\times l}$ be a pattern.
    For $j\in [k]$, define $\Span(j)$ to be the interval of columns between the first and last 1 in 
    row $j$.
    A partition of the rows of $P$ into contiguous intervals $J_0,J_1,\ldots,J_r$
    is a \emph{legal division} 
    if the following properties are satisfied.
\begin{enumerate}
    \item Either $J_0$ is not empty, or at least two of the other intervals are not empty.  
    \item For any two $j,j'\in J_0$, 
    $\Span(j) \cap \Span(j') = \emptyset$.
    \item For every $\ell \in [l]$, all 1s in column $\ell$ lie in $J_0$ and at most one other $J_t$, $t\in [r]$.
    \item If column $\ell$ has a non-zero intersection with $J_t$, $t\in [r]$, then it is \emph{labeled} $t$.
    If it does not intersect any $J_1,\ldots,J_r$ then it is labeled $\perp$.  For every row $j\in J_0$,
    the non-$\perp$ labels of columns in $\Span(j)$ 
    are strictly increasing, regardless of 
    whether the entries here are 0 or 1.
\end{enumerate}
\end{definition}

\begin{theorem}
    If $P$ has no legal division then $P\nprec A^{\mathrm{anti}}$ and $\ex(n,P)=\Omega(n\log n/\log\log n)$.
\end{theorem}

\begin{proof}
    Suppose $P\prec A_{k,r}^{\text{anti}}$.  
    Partition the rows of $A^{\text{anti}}$ that participate in the embedding of $P$ into $J_0,\ldots,J_r$ according to their first digit, and assume without loss of generality that it is not the case that all these rows lie in a 
    single $J_t$, $t\in [r]$. 
    For every $j\in J_0$, the columns 
    containing 1s share a common length-$(k-1)$ suffix.  
    By the antilexicographic ordering, it follows that for any two 
    $j,j' \in J_0$, $\Span(j)$ and $\Span(j')$ do not intersect, for otherwise they would have the same suffix and then $j=j'$.
    Because $A^{\text{anti}}$ only differs from $A$ in the order of its columns, it still holds that every column has its 1s in $J_0$ and at most one other $J_t$, $t\in [r]$.
    If a column intersects $J_t$, $t\in [r]$, 
    then its first digit is $t$, which is also its label.
    For $j\in J_0$, the columns in $\Span(j)$ share a common suffix, so they are sorted according to their first digit/label. 
    In other words, if $P\prec A^{\text{anti}}$ 
    then $P$ has a legal division.
\end{proof}

\begin{corollary}\label{cor:anti-Tardos}
    $S_{12}$ has no legal division, hence
    $S_{12}\nprec A^{\mathrm{anti}}$, hence $\ex(n,S_{12})=\Omega(n\log n/\log\log n)$.    
\[
S_{12} = \patbits{1010/0001/0010/0100/0001}.
\]
\end{corollary}

\begin{proof}
    It cannot be that $J_0 = \emptyset$ since 
    \Cref{def:legal-division}(3) applied to the third and fourth columns would force all rows to 
    be in a single $J_t$, $t\in [r]$.
    It cannot be that the first three rows are in $J_0$,
    as this violates \Cref{def:legal-division}(2): 
    the spans of the first and third rows are not disjoint.
    If the first or first and second rows are in $J_0$,
    then the 1 in the fourth row and the 1 in the third row 
    cannot be assigned distinct labels in increasing order,
    violating \cref{def:legal-division}(4) at row 1.
\end{proof}

\subsection{Upper Bounds}\label{subsect:logn-over-loglogn-ub}

We give a general method for proving 
$O(n\log n/\log\log n)$ bounds on $\ex(n,S)$ if 
$S$ satisfies a few properties. 
For full generality,
\Cref{thm:logn-over-loglogn-decomposition} considers rectangular matrices and \emph{sets} of forbidden patterns.
To keep the statement of \Cref{thm:logn-over-loglogn-decomposition} simple we only consider \emph{superadditive} $\calS$.

\begin{lemma}[Superadditivity]\label{lem:superadditive}
In general $\ex(n,m,\calS)$ may not satisfy superadditivity. 
Call a pattern \emph{NW-SE decomposable} (resp., \emph{SW-NE decomposable}) if it can be partitioned into nonempty quadrants such that the SW and NE quadrants
(resp., NW and SE quadrants) are zero.
Then if it is the case that every pattern in $\calS$ has nonzero exterior rows and columns, and either every pattern in $\calS$ is not NW-SE decomposable or every pattern in $\calS$ is not SW-NE decomposable, then $\ex(n_1,m_1,\calS) + \ex(n_2,m_2,\calS) \leq \ex(n_1+n_2,m_1+m_2,\calS)$.
When $\cal{S} = \{S\}$ contains a single pattern,
it is always true that
$\ex(n_1,m_1,S) + \ex(n_2,m_2,S) \leq \ex(n_1+n_2,m_1+m_2,S)$.
\end{lemma}

\begin{proof}
    Let $A_1,A_2$ be $n_1\times m_1,n_2\times m_2$ $\calS$-free matrices.  If every pattern in $\calS$ is not NW-SE decomposable, then $\left(\begin{array}{cc}A_1&0\\0&A_2\end{array}\right)$ is also $\calS$-free.  The SW-NE decomposable case is symmetric. The singleton case is proved in \cite{PachTardos06}.
\end{proof}

\begin{theorem}\label{thm:logn-over-loglogn-decomposition}
    Let $\calS$ be a \emph{set}
    of patterns.
    Obtain $\calSf$ from $\calS$ by trimming off the last column of each pattern in $\calS$ whenever it has weight 1.
    $\calSl$ is defined similarly with respect to first columns,
    and $\calSm$ is obtained by trimming off both first and last columns with weight 1.
    Suppose $\ex(n,m,X)$ is superadditive, for $X\in\{\calS,\calSf,\calSl,\calSm\}$.
    Then for any integer $K>1$, there exists an allocation 
    $n = n^* + \sum_{i=1}^K n_i$ such that
    \begin{align*}
    \ex(n,m,\calS) &\leq \sum_{i=1}^K \ex(n_i,m/K,\calS) 
        + \ex(n^*,m,\calSf) + \ex(n^*,m,\calSl) + \ex((K-2)n^*,m,\calSm).
    \intertext{Invoking this bound recursively, $\ex(n,m,\calS)$ is bounded as follows.}
    \ex(n,m,\calS) &\leq n + \ex(n,m\ceil{\log_K m},\calSf) + \ex(n,m\ceil{\log_K m},\calSl) + \ex((K-2)n,m\ceil{\log_K m},\calSm).
    \end{align*}
\end{theorem}

\begin{proof}
    Let $A$ be an $n\times m$ $\calS$-free matrix.  
    Partition $A$ into \emph{slabs} of width $m/K$.
    A row is \emph{local} if its 1s appear in one slab and \emph{global} otherwise.
    The 1s in a global row are partitioned into \emph{first}, \emph{last}, and \emph{middle} depending on whether they are in the row's first nonzero slab, 
    last nonzero slab, or a middle nonzero slab.  Let $n_i$ be the number of rows local to slab $i$, $n^*$ be the number of global rows, and $\nf_i,\nl_i,\nm_i$ be the number of global rows that are first, last, and middle in slab $i$, respectively.  It follows that $\sum_i \nf_i = \sum_i \nl_i = n^*$ and that $\sum_i \nm_i \leq (K-2)n^*$.

    Each of the $\nf_i\times m/K$ matrices of global first 1s is $\calSf$-free since each row contains a 1 in a slab to the right.  Similar observations hold for the $\nl_i\times m/K$ and $\nm_i\times m/K$ matrices of global last and middle 1s.  
    The number of global 1s in $A$ in each category is at most:
    \begin{align*}
        \sum_i \ex(\nf_i,m/K,\calSf) &\leq \ex(n^*,m,\calSf),\\
        \sum_i \ex(\nl_i,m/K,\calSl) &\leq \ex(n^*,m,\calSl),\\
        \sum_i \ex(\nm_i,m/K,\calSm) &\leq \ex((K-2)n^*,m,\calSm).
    \end{align*}
    The number of local 1s is at most $\sum_i \ex(n_i,m/K,\calS)$, which 
    completes the recurrence.  If we apply it recursively with the same value of $K$,
    then all the $n^*$ terms at each level sum to at most $n$, and the $m$ terms sum to $m\ceil{\log_K m}$.
    The rows that are never classified as global are local at the last level of recursion; they contribute at most $n$ 1s. This completes the proof.
\end{proof}

\begin{corollary}\label{cor:logn-over-loglogn}
    Whenever $\ex(n,m,\calSm) = O(n+m)$ is linear and 
    $\ex(n,m,\calSf),\ex(n,m,\calSl)$ are quasilinear, 
    $\ex(n,\calS) = O(n\log n/\log\log n)$.
\end{corollary}

\begin{proof}
    By quasilinear we mean that $\ex(n,m,\calSf) = O(n\beta(n,m) + m)$ 
    or $O(m\beta(m,n) + n)$, where $\beta(n,m) = 2^{O(\alpha^t(n,m))}$ 
    is a function of the inverse-Ackermann function.  
    \Cref{thm:logn-over-loglogn-decomposition} implies that
    $\ex(n,m,\calS)$ is
    \begin{align*}
    O\big(n(K + \beta(n,m\log_K m)) + m\log_K m(1 + \beta(m\log_K m, n))\big).
    \end{align*}
    When $n=m$ we choose $K=\log n/\log\log n$,  
    so $\beta(n,m\log_K m) = 2^{O(\alpha^t(n))} < K$
    and $\beta(m\log_K m,n) = O(1)$.
    Therefore $\ex(n,\calS) = O(n\log n/\log\log n)$.
\end{proof}

\begin{corollary}\label{cor:S:logn-over-loglogn}
    For $k\in [10]$, $\ex(n,S_k) = O(n\log n/\log\log n)$.
\end{corollary}

\begin{proof}
    $S_1$ and $S_2$ satisfy 
    $\ex(n,m,\Sf_i) = O(m\alpha(m,n)+n)$~\cite{FurediH92} and $\ex(n,m,\Sl_i),\ex(n,m,\Sm_i)=O(m+n)$.
    Due to \Cref{thm:Qhat-alpha} and \cite{FurediH92}, 
    $S_3,\ldots,S_{10}$ satisfy 
    $\ex(n,m,\Sf_i),\ex(n,m,\Sl_i) = O(n\alpha(n,m)+m)$
    and $\ex(n,m,\Sm_i) = O(n+m)$.
\end{proof}

\Cref{thm:logn-over-loglogn-decomposition} provides a simpler proof of Tardos's upper bound on $\{R_1,R_1^{\varobar}\}$-free matrices.

\begin{corollary}[{Tardos~\cite[Theorem 3.2]{Tardos05}}]\label{cor:Tardos-simplified}
    $\ex(n,\{R_1,R_1^{\varobar}\}) = O(n\log n/\log\log n)$.
    \[
    R_1 = \patbits{110/101},\qquad
    R_1^{\varobar} = \patbits{011/101}.
    \]
\end{corollary}

\begin{proof}
    Observe that when $\calS = \{R_1,R_1^{\varobar}\}$, $\calSf,\calSl,\calSm$ 
    are all linear,
    so $\ex(n,m,\calS) = O(Kn + m\log_K m)$.  
    For square matrices we choose $K=\log n/\log\log n$.
\end{proof}

\Cref{thm:logn-over-loglogn-decomposition} also implies new upper bounds on other matrices, which in the case of $T_1$ and $T_2$ are \emph{not} known to be tight.

\begin{corollary}
    $\ex(n,T_1), \ex(n,T_2) = O(n(\log n/\log\log n)^2)$.
    \[
    T_1 = \patbits{10101/01110}, \qquad 
    T_2 = \patbits{01101/10110}.
    \]
\end{corollary}

\begin{proof}
    For $S = T_1$ or $T_2$, $\ex(n,m,\Sm) = O(n+m)$ and $\ex(n,m,\Sf),\ex(n,m,\Sl)$
    are both at most $\ex(n,m,R_1)+n$,
    which is $\Theta(\bar{m}\log_{1+\bar{m}/\bar{n}}\bar{m})$, where $\bar{m} = \max\{m,n\}$ and $\bar{n} = \min\{m,n\}$.  The $O(n\log n)$ bounds on $\ex(n,R_1)$ are well known~\cite{Furedi90,FurediH92,Tardos05}.  
    The rectangular bounds are from Pettie~\cite[Corollary A.2]{Pettie10a}.
    Applying \Cref{thm:logn-over-loglogn-decomposition}, we bound $\ex(n,m,S)$ as follows.
    \begin{align*}
    \ex(n,m,S) &\leq \ex(n,m\log_K m,\Sf) + \ex(n,m\log_K m,\Sl) + \ex((K-2)n,m\log_K m,\Sm).
    \intertext{When $m=n$ we set $K=\log n$.  
    Letting $d = \log_K m$, this is bounded by }
        &= O(Kn + m\log_K m\log_{1+d} (m\log_K m)) = O(n(\log n/\log\log n)^2).
    \end{align*}
\end{proof}

\section{Conclusion}\label{sec:conclusion}

We have fully classified all 
weight-5 patterns in terms of their linear/nonlinear status
and made significant progress on classifying weight-5 and weight-6 patterns according to their extremal functions, asymptotically.
Some notable results are several new infinite classes
of linear patterns (captured by Theorems~\ref{thm:three-ones}, \ref{thm:two-sided-extension}, \ref{thm:one-sided-extension}, and \ref{thm:hk-minusone-linear}), 
a new class of patterns with extremal function
$\Theta(n\log n/\log\log n)$ (\Cref{thm:logn-over-loglogn-decomposition}, 
and \Cref{cor:S-one-through-ten,cor:S:logn-over-loglogn}) and a $\Theta(n\alpha(n))$ bound on a nonlight pattern (\Cref{thm:Qhat-alpha}).
To complete the classification at weight-5, 
we need to understand the following patterns, 
which we group into five echelons.

\bigskip 
The following two light patterns  are known to be $\Omega(n\alpha(n))$ and $O(n\alpha^2(n))$~\cite{FurediH92,Pettie-GenDS11}.
We think the lower bound is more likely to be tight.  
On the other hand, there are matrix constructions with 
weight $O(n\alpha^t(n))$ that avoid some structured patterns.
See \cite[\S 8]{Pettie15-SIDMA}.
\[
\patword{31323}, \qquad \patword{41323}.
\]
There are three light patterns known to be 
$\Omega(n\alpha(n))$ and $O(n2^{\alpha(n)})$~\cite{FurediH92,Pettie-GenDS11,ChalermsookPY24,Pettie15-SIDMA}.\footnote{The patterns \emph{not} on this list are due to
recent $\Theta(n2^{\alpha(n)})$ bounds~\cite{ChalermsookPY24},
new $\Theta(n\alpha(n))$ bounds 
(\Cref{thm:Qtilde-alpha}),
new linear bounds on $Q_{41423},Q_{41432},Q_{42314}$, 
and patterns certified $\Theta(n\alpha(n))$ by 
\Cref{thm:one-sided-extension} applied to a known $\Theta(n\alpha(n))$ pattern.  For example, starting
with a pattern of Pettie~\cite{Pettie15-SIDMA},
\[
\patword{21221} 
\xrightarrow{\,\mbox{\small trivial}\,}
\patword{21321}
\xrightarrow{\,\mbox{\small Thm.~\ref{thm:one-sided-extension}+trivial}\,}
\patword{31421}.
\]}
An $O(n\alpha(n))$ upper bound on the first would 
resolve the other two.
\[
\begin{array}{c@{\quad}c@{\quad}c@{\quad}c@{\quad}c@{\quad}c}
\patword{31321}, \qquad
\patword{42431}, \qquad
\patword{41321}.
\end{array}
\]
The following three patterns are \emph{not} light, 
and are known to be $\Omega(n\alpha(n))$
as they all contain $Q_3'$~\cite{FurediH92}.
The first is known to be $O(n\log n)$~\cite{Tardos05,PachTardos06},
while \Cref{thm:logn-over-loglogn-decomposition} and 
\Cref{cor:logn-over-loglogn} imply 
the other two are $O(n\log n/\log\log n)$.
\[
\begin{array}{c@{\quad}c@{\quad}c@{\quad}c}
\patbits{0001/0100/0001/1010},\qquad
\patbits{1000/0010/0101/0010},\qquad
\patbits{0100/0010/0100/1001}.
\end{array}
\]
The following patterns all contain $R_1,R_2$, 
or $R_3$, and are therefore $\Omega(n\log n)$.
They are all $O(n(\log n/\log\log n)^2)$ due to~\cite{PettieT25}.
\[
\patbits{1010/1100/0001},\quad
\patbits{1010/0101/1000},\quad
\patbits{1010/0100/1001},\quad
\patbits{1010/0100/1000/0001},\quad
\patbits{1010/0100/0001/1000}.
\]
Finally, there are three patterns known to be 
$\Omega(n\log n)$ and $O(n\log^2 n)$~\cite{PachTardos06}.
\[
\patbits{110/001/101},\quad\quad
\patbits{110/101/010},\quad\quad
\patbits{110/010/101}.
\]

\section*{Declaration of AI Use}
We used Codex with GPT-6 Astra Ultra to enumerate small 
weight patterns with stated properties, 
to assist with finding proofs, 
to generate \Cref{fig:large-LR-diagram,fig:composition-patterns}, 
and for mathematical copyediting.  
The paper \emph{per se} is entirely 
written by the authors.

\bibliographystyle{alpha}
\bibliography{references}

\appendix

\section{New Quasilinear Bounds}\label{sect:new-quasilinear-bounds} 

\subsection{A Quasilinear Bound on a Nonlight Pattern}\label{subsect:quasilinear-nonlight-pattern}

All light patterns are quasilinear~\cite{Klazar92,Nivasch10,Pettie15-SIDMA},
but not all quasilinear patterns are light.
One can construct somewhat artificial 
nonlight patterns that are quasilinear using Keszegh's~\cite{Keszegh09} gluing operation; see~\cite[Footnote 2]{PettieT24}.
\Cref{thm:Qhat-alpha}
establishes a sharp $\Theta(n\alpha(n,m) + m)$
bound on the ``genuinely'' nonlight pattern $\hat{Q}_3$. The lower bound follows from~\cite{HS86,FurediH92}.
\[
\hat{Q}_3 = \patbits{1010/0100/0101}.
\]

\Qhat*

\begin{proof}
    Let $A$ be an $n\times m$, $\hat{Q}_3$-free 0--1 matrix.
    Partition $A$ into \emph{slabs} of $B$ consecutive columns.  Every row is \emph{local} if its 1s appear in a single slab and \emph{global} otherwise.  
    The 1s of a global row are partitioned into \emph{first}, \emph{last}, and \emph{middle} occurrences, according to whether they appear in the first non-zero slab, the last non-zero slab, or any other non-zero slab intersecting the row.
    Let $n_i$ be the number of rows local to slab $i$
    and $n^*$ be the number of global rows. Define $\nf_i,\nm_i,\nl_i$ to be the number of rows that have a global first, middle, and last 1 in slab $i$, respectively.  Define the subpatterns 
    $\Uf,\Ul,\Um$ as:
    \[
    \Uf = \patbits{101/010/010},
    \qquad
    \Ul = \patbits{010/100/101},
    \qquad
    \Um = \patbits{01/10/10}.
    \]
    In each slab, the submatrix of first 1s is $\Uf$-free and the submatrix of last 1s is $\Ul$-free.
    The number of these 1s is therefore
    $\sum_{i} (\ex(\nf_i,B,\Uf) + \ex(\nl_i,B,\Ul)) = O(n^* + m)$.
    The submatrix of middle 1s in each slab is $\Um$-free.
    Remove the first and last global middle 1 in each column,
    and let $\bar{A}$ be the matrix of the remaining middle 1s.
    It follows that in each slab $i$, no row of $\bar{A}$ 
    has two 1s.  Define $\bar{A}'$ to be the $n^* \times m/B$ 
    matrix obtained by contracting each slab to a single column and retaining only global rows. It follows from the previous observation that $|\bar{A}| = |\bar{A}'|$.
    
    In general, the class of matrices avoiding a nonlight pattern is not closed under contracting columns.  Nonetheless, we claim $\bar{A}'$ is $\hat{Q}_3$-free.  
    Consider an occurrence of $\hat{Q}_3$ in $\bar{A}'$.
    \[
    \left(\begin{array}{c|c|cc}
        \bu&&\bu&\\
        &\bu&&\\
        &\bu&&\bu
    \end{array}\right)
    \]
    The 1s in the second column correspond to two middle 1s in a slab of $\bar{A}$, say 
    $\bar{A}(i_0,j_0) = \bar{A}(i_1,j_1) = 1$, where $i_0 < i_1$. 
    It cannot be that $j_0 = j_1$ for this implies that 
    $\hat{Q}_3$ appears in $A$. 
    It cannot be that $j_0 > j_1$ for this implies that the slab
    was not $\Um$-free. See the diagram on the left below; the $\circ$ is the last 1 removed from column $j_1$.  It must be that $j_0 < j_1$, as in the diagram on the right below.  Let $A(i_2,j_1)=1$ be the top 1 in column $j_1$ that was removed.
    Then $i_2 \in [i_0,i_1)$, for otherwise the slab would not be $\Um$-free.  However, this implies an occurrence of $\hat{Q}_3$ in the original matrix $A$, substituting row $i_2$ for $i_0$.
    \bigskip 
    \[
    \begin{array}{r}
    i_0\\i_1\\\ 
    \end{array}\hcm[-.2]
    \left(\begin{array}{cc}
    \zero{\rb{5}{$j_1$}} & \zero{\rb{5}{$\:j_0$}}\bu\\
    \bu&\\
    \circ&
    \end{array}\right),
    \qquad\qquad\qquad
    \begin{array}{r}
    i_0\\i_2\\i_1 
    \end{array}\hcm[-.2]
    \left(\begin{array}{cc}
    \zero{\rb{5}{$j_0$}}\bu & \zero{\rb{5}{$\:j_1$}}\\
    &\circ\\
    &\bu
    \end{array}\right).
    \]
    Since $A$ was an arbitrary $\hat{Q}_3$-free matrix, 
    it follows that
    \[
    \ex(n,m,\hat{Q}_3) \leq \sum_i \ex(n_i,B,\hat{Q}_3) + \ex(n^*,m/B,\hat{Q}_3) + O(n^* + m).
    \]
    Whenever a function satisfies this recurrence, it can be proved that $\ex(n,m,\hat{Q}_3) = O(n\alpha(n,m) +m)$. 
    See~\cite{HS86,Pettie-GenDS11} for example.
\end{proof}

\subsection{Sharp Bounds on a Light Pattern}\label{subsect:quasilinear-light-pattern}

We prove sharp $\Theta(n\alpha(n))$ bounds for $\tilde{Q}_3$-free matrices.  The lower bound follows from~\cite{HS86,FurediH92}.
\[
\tilde{Q}_3 = \patword{31231}.
\]

\Qtilde*

\begin{proof}
    Let $A_0$ be an $n\times m$ $\tilde{Q}_3$-free matrix. Partition $A_0$ into slabs of width $B$ and classify its rows as local or global,
    where $n^*$ is the number of global rows.
    Obtain $A$ by removing the last 1 in each global row of $A_0$, so $|A|=|A_0|-n^*$.
    Partition the global 1s of $A$ into first, middle, and last, and define $n^*,n_i,\nf_i,\nl_i,\nm_i$ 
    as in the proof of \Cref{thm:Qhat-alpha}, 
    with one exception.  Note that a row $r$ is classified as global with respect to $A_0$, but removing the last 1 from row $r$ could cause its remaining 1s to lie in a single slab, say slab $i$. In this case $r$ would be counted in $\nf_i$ but not $\nl_i$.
    Since trimming off the first or last column of $\tilde{Q}_3$ makes the extremal function linear,
    it follows that the contribution of global first and global last 1s is $O(n^* + m)$.  
    The number of local 1s is at most $\sum_i \ex(n_i,B,\tilde{Q}_3)$.
    We now analyze $\bar{A}$, the matrix of global middle 1s.

    Each slab of $\bar{A}$ is $U_3$-free, for otherwise $A$ would contain $\tilde{Q}_3$.
    \[
    U_3 = \patword{123}.
    \]
    Classify the middle 1s in each slab of $\bar{A}$ as type-0 if the slab contains no middle 1 that is strictly to the NW of it, 
    and type-1 otherwise.  It follows that the submatrix of type-0 1s in the slab is $U_2$-free, and the submatrix of type-1 1s in the slab is $U_2$-free.
    \[
    U_2 = \patword{12}, \qquad U_2^{\text{rot}} = \patword{21}.
    \]
    The type-0 matrix is $U_2$-free by definition.  The type-1 matrix cannot contain a $U_2$, for otherwise the slab of $\bar{A}$ would then contain a forbidden $U_3$.
    Delete the top type-0 and type-1 1s from each column of $\bar{A}$,
    which leaves at most one type-0 and type-1 1 per row.
    The matrix $\bar{A}'$ is formed from the remaining 1s by contracting each slab into \emph{two} adjacent columns: the first retains all type-0 1s and the second all type-1 1s.  An example is given below, with type-0 and type-1 marked as $\circ$ and $\bu$.
    \[
    \left(\begin{array}{cccccccc}
    &&&&&&\circ&\\
    &&&&&&&\bu\\
    &&&&\circ&&&\\
    &&&&&&\bu&\\
    &&&&&\bu&&\\
    &&\circ&&&&&\\
    &&&\bu&&&&\\
    \circ&&&&&&&\\
    &\bu&&&&&
    \end{array}\right)
    \qquad
    \longrightarrow
    \qquad
    \left(\begin{array}{c|c}
        \circ &\\
        &\bu\\
     \circ &\\
        &\bu\\
        &\bu\\
                \circ &\\
        &\bu\\
                \circ &\\
        &\bu
    \end{array}\right)
    \]
    
    Note that $\bar{A}'$ is \underline{\emph{not}} a column-minor of 
    $\bar{A}$, so we cannot immediately conclude that it inherits all the light forbidden patterns of $\bar{A}$.  
    Nonetheless, we claim that $\bar{A}'$ is $\tilde{Q}_3$-free.
    
    Let us first examine the ways that two 1s from the slab 
    of $\bar{A}$ can be arranged in $\bar{A}'$ when they 
    have opposite types.
    Suppose rows $x$ and $y$ of $\bar{A}'$ retain type-0 and type-1 1s from a common slab in $\bar{A}$, in columns $j_x$ and $j_y$, respectively.  If $x$ is below $y$ then $j_x < j_y$, i.e.,
    their order in $\bar{A}$ agrees with their order in $\bar{A}'$.
    If $j_y \leq j_x$, 
    then $\bar{A}(x,j_x)=1$ would have been 
    classified as type-1.
    \[
    \begin{array}{r}
    \ \\
    \ \\
    \ \\
    y\\
    \ \\
    \ \\
    x\\\ 
    \end{array}
    \left(\begin{array}{c|c}
    &\\
    &\\
    &\\
    &\bu\\
    &\\
    &\\
    \circ&\\
    \multicolumn{1}{r|}{\zero{\rb{-5}{\hcm[-.9]type-0}}} 
    & \multicolumn{1}{l}{\zero{\rb{-5}{\hcm[.9]type-1}}}\\
    \end{array}\right)
    \qquad\longleftarrow\qquad
    \left(\begin{array}{cccccc}
    \multicolumn{2}{c|}{\hcm[1.5]} & j_x & \multicolumn{1}{c|}{\hcm[1.5]} & j_y&\\
    \multicolumn{2}{c|}{\emptyset} & \multicolumn{2}{c|}{\neq\emptyset} &&\\
    \multicolumn{2}{c|}{} & \multicolumn{2}{c|}{} &&\\\cline{1-4} 
    \multicolumn{2}{c|}{} & & & \bu &\\
    \multicolumn{2}{c|}{\emptyset} & & & &\\
    \multicolumn{2}{c|}{} & & & &\\\cline{1-2}
    &&\circ& & &\\
    &&&&&\\
    \end{array}\right)
    \]

\bigskip
    When $x$ is above $y$  then we can have either $j_x<j_y$
    (agreeing with $\bar{A}$) or $j_y \leq j_x$, but in the latter case we can deduce the existence of a third point.
    \[
    \begin{array}{r}
    \ \\
    \ \\
    \ \\
    x\\
    \ \\
    \ \\
    y\\\ 
    \end{array}
    \left(\begin{array}{c|c}
    &\\
    &\\
    &\\
    \circ&\\
    &\\
    &\\
    &\bu\\
    \multicolumn{1}{r|}{\zero{\rb{-5}{\hcm[-.9]type-0}}} 
    & \multicolumn{1}{l}{\zero{\rb{-5}{\hcm[.9]type-1}}}\\
    \end{array}\right)
    \qquad\longleftarrow\qquad
    \left(\begin{array}{cccccc}
    \multicolumn{2}{c|}{\hcm[1.5]} & j_y & \multicolumn{1}{c|}{\hcm[1.5]} & j_x&\\
    \multicolumn{2}{c|}{\emptyset} & \multicolumn{2}{c|}{\emptyset} &&\\
    \multicolumn{2}{c|}{} & \multicolumn{2}{c|}{} &&\\\cline{1-4} 
    \multicolumn{2}{c|}{} & & & \circ &\\
    \multicolumn{2}{c|}{\circ} & & & &\\
    \multicolumn{2}{c|}{} & & & &\\\cline{1-2}
    &&\bu& & &\\
    &&&&&\\
    \end{array}\right)
    \]

\bigskip 
In other words, whenever we see a $U_2^{\text{rot}}$ in adjacent columns in $\bar{A}'$, it corresponds to a $U_2^{\text{rot}}$ in the corresponding slab of $\bar{A}$. 
Whenever we see a $U_2$ in $\bar{A}'$
on rows $x<y$, it corresponds to a $U_2$ in $\bar{A}$
on rows $z,y$ where $z\in [x,y)$.

Now suppose that $\bar{A}'$ contains an occurrence of $\tilde{Q}_3$ 
on rows $c < b < a$ and columns generated by 
slabs $i_1 \leq \cdots \leq i_5$, adjacent pairs of which may be equal.
\medskip 

\[
\tilde{Q}_3 = 
\begin{array}{l}
    c\\b\\a\\
\end{array}
\hcm[-.2]
\left(\begin{array}{ccccc}
\zero{\rb{4}{$i_1$}}&\zero{\rb{4}{$i_2$}}\bu&\zero{\rb{4}{$i_3$}}&\zero{\rb{4}{$i_4$}}&\zero{\rb{4}{$i_5$}}\bu\\
    &&\bu&&\\
\bu&&&\bu&
\end{array}\right)
\]
If $i_1,\ldots,i_5$ are distinct, or if the only equalities are $i_1=i_2$ or $i_4=i_5$, then this occurrence of $\tilde{Q}_3$ also exists in $\bar{A}$.  The interesting cases are if $i_2=i_3$ or $i_3=i_4$.  When $i_3=i_4$, the 1 in row $b$ must be type-0 while the 1 in row $a$ is type-1.  In slab $i_3$, there is either a copy of $U_2$ on rows $b,a$ or a copy on rows $d,a$, for some $d \in (b,a)$.  
In the latter case, there is an occurrence of $\tilde{Q}_3$ 
on rows $a,d,c$.
\[
\tilde{Q}_3 = 
\begin{array}{l}
    c\\b\\d\\a\\
\end{array}
\hcm[-.2]
\left(\begin{array}{cc|ccc|c}
\zero{\rb{4}{$i_1$}}&\zero{\rb{4}{$i_2$}}\bu&&\zero{\rb{4}{$i_3$}}&&\zero{\rb{4}{$i_5$}}\bu\\
    &&&&\circ&\\
    &&\circ&&&\\
\bu&&&\bu&
\end{array}\right)
\]
If $i_2=i_3$ then slab $i_2$ either contains a copy of $U_2$ on rows $b,c$ (which implies an occurrence of $\tilde{Q}_3$ in $\bar{A}$)
or a copy of $U_2$ on rows $b,d$ for some $d\in (c,b)$.  
This is a global middle 1 in row $d$, 
which implies there are at least two 1s in row $d$ to the right: one in the \emph{last} slab for row $d$, and another from $A_0$ not appearing in $A$.  Thus, there is a 1 in row $d$ that is not in the same column as the 1 at the intersection of row $a$ and slab $i_4$ of $\bar{A}$.
There are two cases, depending on whether the 1 in row $d$ 
appears strictly before or strictly after the 1 in row $a$.
$A_0$ contains $\tilde{Q}_3$ on rows $a,d,c$ in the first case
and rows $a,b,d$ in the second, as depicted below.
\end{proof}

\[
\begin{array}{l}
    c\\d\\b\\a\\
\end{array}
\hcm[-.2]
\left(\begin{array}{c|ccc|ccc}
&&\zero{\rb{4}{$i_2$}}&\circ&&&\bu\\
&\circ&&&\bu&&\\
&&\bu&&&&\\
\bu&&&&&\bu&
\end{array}\right),
\qquad\qquad
\begin{array}{l}
    c\\d\\b\\a\\
\end{array}
\hcm[-.2]
\left(\begin{array}{c|ccc|ccc}
&&\zero{\rb{4}{$i_2$}}&\circ&&&\bu\\
&\circ&&&&\bu&\\
&&\bu&&&&\\
\bu&&&&\bu&
\end{array}\right).
\]

\medskip
Putting it all together, $\ex(n,m,\tilde{Q}_3)$ obeys
the following recurrence:
\[
\ex(n,m,\tilde{Q}_3) \leq \sum_i \ex(n_i,B,\tilde{Q}_3) + \ex(n^*,2m/B,\tilde{Q}_3) + O(n^* + m).
\]
The only difference between this recurrence and the standard one
(see \Cref{thm:Qhat-alpha})
is that $\bar{A}'$ contains $2m/B$ columns 
instead of $m/B$.  
This does not affect the $O(n\alpha(n,m)+m)$ bound 
asymptotically; see~\cite{Pettie-GenDS11,Pettie15-SIDMA}.
\end{document}

%% file: patbits.tex
\newif\ifpatfirst
\newcount\patmax
\newcount\patrow
\newcount\patdigit
\newtoks\patbody

\def\patappend#1{%
  \patbody=\expandafter{\the\patbody #1}%
}

\def\patmatrix{%
  \mathord{%
    \left(%
      \vcenter{%
        \offinterlineskip
        \halign{%
          \strut\hfil$##$\hfil
          &&%
          \quad\strut\hfil$##$\hfil
          \cr
          \the\patbody
          \cr
        }%
      }%
    \right)%
  }%
}

\def\patmatrixcompact{%
  \mathord{%
    \left(%
      \vcenter{%
        \offinterlineskip
        \halign{%
          \strut\hfil$##$\hfil
          &&%
          \strut\hfil$##$\hfil 
          \cr
          \the\patbody
          \cr
        }%
      }%
    \right)%
  }%
}

\def\patbitsscan#1{%
  \ifx#1\relax
    \let\patnext\relax
  \else
    \ifx#1/%
      \patappend{\cr}%
      \patfirsttrue
    \else
      \ifpatfirst
        \patfirstfalse
      \else
        \patappend{&}%
      \fi
      \ifx#11%
        \patappend{\bu}%
      \else
        \patappend{\phantom{\bu}}%
      \fi
    \fi
    \let\patnext\patbitsscan
  \fi
  \patnext
}

\def\patbits#1{%
  \begingroup
  \patbody={}%
  \patfirsttrue
  \patbitsscan#1\relax
  \patmatrix
  \endgroup
}

\def\patfindmax#1{%
  \ifx#1\relax
    \let\patnext\relax
  \else
    \patdigit=#1\relax
    \ifnum\patdigit>\patmax
      \patmax=\patdigit
    \fi
    \let\patnext\patfindmax
  \fi
  \patnext
}

\def\patwordscan#1{%
  \ifx#1\relax
    \let\patnext\relax
  \else
    \ifpatfirst
      \patfirstfalse
    \else
      \patappend{&}%
    \fi
    \patdigit=#1\relax
    \ifnum\patdigit=\patrow
      \patappend{\bu}%
    \else
      \patappend{\phantom{\bu}}%
    \fi
    \let\patnext\patwordscan
  \fi
  \patnext
}

\def\patwordrows#1{%
  \ifnum\patrow>\patmax
  \else
    \ifnum\patrow>1
      \patappend{\cr}%
    \fi
    \patfirsttrue
    \patwordscan#1\relax
    \advance\patrow by 1
    \patwordrows{#1}%
  \fi
}

\def\patword#1{%
  \begingroup
  \patbody={}%
  \patmax=0
  \patfindmax#1\relax
  \patrow=1
  \patwordrows{#1}%
  \patmatrix
  \endgroup
}

\def\patwordc#1{%
  \begingroup
  \patbody={}%
  \patmax=0
  \patfindmax#1\relax
  \patrow=1
  \patwordrows{#1}%
  \patmatrixcompact
  \endgroup
}

%% file: large-diagram.tex
\newcommand{\zerobox}[2][0pt]{%
  \smash{\makebox[0pt][l]{%
    \hspace*{#1}\makebox[0pt][c]{$#2$}%
  }}%
}

\setlength{\arraycolsep}{1.1pt}
\renewcommand{\arraystretch}{0.6}
\setlength{\arrayrulewidth}{0.35pt}
\setlength{\dashlinedash}{2pt}
\setlength{\dashlinegap}{1.5pt}
\ADLnullwide

\setbox\strutbox=\hbox{\vrule height10.16667pt depth1.83333pt width0pt}%
\(
\begin{array}{*{56}{c}}

\cline{1-12}
  \multicolumn{1}{|c}{\hspace*{.5em}} & \hspace*{.5em} & \hspace*{.5em} & \hspace*{.5em} & \hspace*{.5em} & \bu & \bu & &
   & \hspace*{.5em} & \hspace*{.5em} & \multicolumn{1}{c|}{\hspace*{.5em}} & \hspace*{.5em} & \hspace*{.5em} & \hspace*{.5em} & \hspace*{.5em} &
  \hspace*{.5em} & \hspace*{.5em} & \hspace*{.5em} & \hspace*{.5em} & \hspace*{.5em} & \hspace*{.5em} & \hspace*{.5em} & \hspace*{.5em} &
  \hspace*{.5em} & \hspace*{.5em} & \hspace*{.5em} & \hspace*{.5em} & \hspace*{.5em} & \hspace*{.5em} & \hspace*{.5em} & \hspace*{.5em} &
  \hspace*{.5em} & \hspace*{.5em} & \hspace*{.5em} & \hspace*{.5em} & \hspace*{.5em} & \hspace*{.5em} & \hspace*{.5em} & \hspace*{.5em} &
  \hspace*{.5em} & \hspace*{.5em} & \hspace*{.5em} & \hspace*{.5em} & \hspace*{.5em} & \hspace*{.5em} & \hspace*{.5em} & \hspace*{.5em} &
  \hspace*{.5em} & \hspace*{.5em} & \hspace*{.5em} & \hspace*{.5em} & \hspace*{.5em} & \hspace*{.5em} & \hspace*{.5em} & \hspace*{.5em} \\
  \multicolumn{1}{|c}{} &  &  &  &  &  &  &  &
   &  &  & \multicolumn{1}{c|}{} &  &  &  &  &
   &  &  &  &  &  &  &  &
   &  &  &  &  &  &  &  &
   &  &  &  &  &  &  &  &
   &  &  &  &  &  &  &  &
   &  &  &  &  &  &  &  \\
  \multicolumn{1}{|c}{} &  &  &  &  &  &  &  &
   &  &  & \multicolumn{1}{c|}{} &  &  &  &  &
   &  &  &  &  &  &  &  &
   &  &  &  &  &  &  &  &
   &  &  &  &  &  &  &  &
   &  &  &  &  &  &  &  &
   &  &  &  &  &  &  &  \\
  \multicolumn{1}{|c}{} &  &  &  &  &  &  & \zerobox{A'_{\mathrm{N}}} &
   &  &  & \multicolumn{1}{c|}{} &  &  &  &  &
   &  &  &  &  &  &  &  &
   &  &  &  &  &  &  &  &
   &  &  &  &  &  &  &  &
   &  &  &  &  &  &  &  &
   &  &  &  &  &  &  &  \\
  \multicolumn{1}{|c}{} &  &  &  &  &  &  &  &
   &  &  & \multicolumn{1}{c|}{} &  &  &  &  &
   &  &  &  &  &  &  &  &
   &  &  &  &  &  &  &  &
   &  &  &  &  &  &  &  &
   &  &  &  &  &  &  &  &
   &  &  &  &  &  &  &  \\
  \multicolumn{1}{|c}{} &  &  &  &  &  &  &  &
   &  &  & \multicolumn{1}{c|}{} &  &  &  &  &
  \multicolumn{1}{c:}{} &  &  &  &  &  &  &  &
   &  &  &  &  &  &  &  &
   &  &  &  &  &  &  &  &
   &  &  &  &  &  &  &  &
   &  &  &  &  &  &  &  \\
  \multicolumn{1}{|c}{} &  &  &  & \bu &  &  &  &
   &  &  & \multicolumn{1}{c|}{} &  &  &  &  &
  \multicolumn{1}{c:}{} &  &  &  &  &  &  &  &
   &  &  &  &  &  &  &  &
   &  &  &  &  &  &  &  &
   &  &  &  &  &  &  &  &
   &  &  &  &  &  &  &  \\
\cline{1-24}
   &  &  &  &  &  &  &  &
   &  &  &  & \multicolumn{1}{|c}{} &  &  &  &
  \multicolumn{1}{c:}{} & \bu &  \bu &  & &  &  & \multicolumn{1}{c|}{} &
   &  &  &  &  &  &  &  &
   &  &  &  &  &  &  &  &
   &  &  &  &  &  &  &  &
   &  &  &  &  &  &  &  \\
   &  &  &  &  &  &  &  &
   &  &  &  & \multicolumn{1}{|c}{} &  &  &  &
  \multicolumn{1}{c:}{} &  &  &  &  &  &  & \multicolumn{1}{c|}{} &
   &  &  &  &  &  &  &  &
   &  &  &  &  &  &  &  &
   &  &  &  &  &  &  &  &
   &  &  &  &  &  &  &  \\
   &  &  &  &  &  &  &  &
   &  &  &  & \multicolumn{1}{|c}{} &  &  &  &
  \multicolumn{1}{c:}{} &  &  &  &  &  &  & \multicolumn{1}{c|}{} &
   &  &  &  &  &  &  &  &
   &  &  &  &  &  &  &  &
   &  &  &  &  &  &  &  &
   &  &  &  &  &  &  &  \\
   &  &  &  &  &  &  &  &
   &  &  &  & \multicolumn{1}{|c}{} &  &  &  &
  \multicolumn{1}{c:}{} &  &  &  &  &  &  & \multicolumn{1}{c|}{} &
   &  &  &  &  &  &  &  &
   &  &  &  &  &  &  &  &
   &  &  &  &  &  &  &  &
   &  &  &  &  &  &  &  \\
   &  &  &  &  &  &  &  &
   &  &  &  & \multicolumn{1}{|c}{} &  &  &  &
  \multicolumn{1}{c:}{} &  &  & \zerobox{A_{\mathrm{N}}} &  &  &  & \multicolumn{1}{c|}{} &
   &  &  &  &  &  &  &  &
   &  &  &  &  &  &  &  &
   &  &  &  &  &  &  &  &
   &  &  &  &  &  &  &  \\
   &  &  &  &  &  &  &  &
   &  &  &  & \multicolumn{1}{|c}{} &  &  &  &
  \multicolumn{1}{c:}{} &  &  &  &  &  &  & \multicolumn{1}{c|}{} &
   &  &  &  &  &  &  &  &
   &  &  &  &  &  &  &  &
   &  &  &  &  &  &  &  &
   &  &  &  &  &  &  &  \\
   &  &  &  &  &  &  &  &
   &  &  &  & \multicolumn{1}{|c}{} &  &  &  &
  \multicolumn{1}{c:}{\bu} &  &  &  &  &  &  & \multicolumn{1}{c|}{} &
   &  &  &  &  &  &  & \zerobox{z} &
   &  &  &  &  &  &  &  &
   &  &  &  &  &  &  &  &
   &  &  &  &  &  &  &  \\
\cline{13-24}
   &  &  &  &  &  &  &  &
   &  &  &  &  &  &  &  &
  \multicolumn{1}{c:}{} &  &  &  &  &  &  &  &
   &  &  &  &  &  & \bu &  &
   &  &  &  &  &  &  &  &
   &  &  &  &  &  &  &  &
   &  &  &  &  &  &  &  \\
\cline{13-24}
   &  &  &  &  &  &  &  &
   &  &  &  & \multicolumn{1}{|c}{} &  &  &  &
  \multicolumn{1}{c:}{} &  &  &  &  &  &  & \multicolumn{1}{c|}{} &
   &  &  &  &  &  &  &  &
   &  &  &  &  &  &  &  &
   &  &  &  &  &  &  &  &
   &  &  &  &  &  &  &  \\
   &  &  &  &  &  &  &  &
   &  &  &  & \multicolumn{1}{|c}{} &  &  &  &
  \multicolumn{1}{c:}{} &  &  &  &  &  &  & \multicolumn{1}{c|}{} &
   &  &  &  &  &  &  &  &
   &  &  &  &  &  &  &  &
   &  &  &  &  &  &  &  &
   &  &  &  &  &  &  &  \\
   &  &  &  &  &  &  &  &
   &  &  &  & \multicolumn{1}{|c}{} &  &  &  &
  \multicolumn{1}{c:}{} &  &  & \zerobox{A_{\mathrm{S}}} &  &  &  & \multicolumn{1}{c|}{} &
   &  &  &  &  &  &  &  &
   &  &  &  &  &  &  &  &
   &  &  &  &  &  &  &  &
   &  &  &  &  &  &  &  \\
   &  &  &  &  &  &  &  &
   &  &  &  & \multicolumn{1}{|c}{} &  &  &  &
  \multicolumn{1}{c:}{} &  &  &  &  &  &  & \multicolumn{1}{c|}{} &
   &  &  &  &  &  &  &  &
   &  &  &  &  &  &  &  &
   &  &  &  &  &  &  &  &
   &  &  &  &  &  &  &  \\
   &  &  &  &  &  &  &  &
   &  &  &  & \multicolumn{1}{|c}{} &  &  &  &
  \multicolumn{1}{c:}{} &  &  &  &  &  &  & \multicolumn{1}{c|}{} &
   &  &  &  &  &  &  &  &
   &  &  &  &  &  &  &  &
   &  &  &  &  &  &  &  &
   &  &  &  &  &  &  &  \\
\cline{13-30}\cline{32-36}
   &  &  &  &  &  &  &  &
   &  &  &  &  &  &  &  &
  \multicolumn{1}{c:}{} &  &  &  &  &  &  &  &
  \multicolumn{1}{|c}{} &  &  &  &  & \multicolumn{1}{c|}{} &  & \multicolumn{1}{|c}{} &
   &  &  & \multicolumn{1}{c|}{} &  &  &  &  &
   &  &  &  &  &  &  &  &
   &  &  &  &  &  &  &  \\
   &  &  &  &  &  &  &  &
   &  &  &  &  &  &  &  &
  \multicolumn{1}{c:}{} &  &  &  &  &  &  &  &
  \multicolumn{1}{|c}{} &  &  &  &  & \multicolumn{1}{c|}{} &  & \multicolumn{1}{|c}{} &
   &  &  & \multicolumn{1}{c|}{} &  &  &  &  &
   &  &  &  &  &  &  &  &
   &  &  &  &  &  &  &  \\
   &  &  &  &  &  &  &  &
   &  &  &  &  &  &  &  &
  \multicolumn{1}{c:}{} &  &  &  &  &  &  &  &
  \multicolumn{1}{|c}{} &  &  & \zerobox{\bar{B}_{\mathrm{NW}}} &  & \multicolumn{1}{c|}{} &  & \multicolumn{1}{|c}{} &
   & \zerobox{\bar{B}_{\mathrm{NE}}} &  & \multicolumn{1}{c|}{} &  &  &  &  &
   &  &  &  &  &  &  &  &
   &  &  &  &  &  &  &  \\
   &  &  &  &  &  &  &  &
   &  &  &  &  &  &  &  &
  \multicolumn{1}{c:}{} &  &  &  &  &  &  &  &
  \multicolumn{1}{|c}{} &  &  &  &  & \multicolumn{1}{c|}{} &  & \multicolumn{1}{|c}{} &
   &  &  & \multicolumn{1}{c|}{} &  \zerobox{e} &  &  &  &
   &  &  &  &  &  &  &  &
   &   &  &  &  &  &  &  \\
\cline{15-21}\cline{36-39}
   &  &  &  &  &  &  &  &
   &  &  &  &  &  & \multicolumn{1}{|c}{} &  &
  \multicolumn{1}{c:}{} &  &  &  & \multicolumn{1}{c|}{} &  &  &  &
  \multicolumn{1}{|c}{} &  &  &  &  & \multicolumn{1}{c|}{} &  & \multicolumn{1}{|c}{} &
   &  &  & \multicolumn{1}{|c|}{\circ} &  &  & \multicolumn{1}{c|}{} &  &
   &  &  &  &  &  &  &  &
   &  &  &  &  &  &  &  \\
%
\cline{25-30}\cline{32-36}
   &  &  &  &  &  &  &  &
   &  &  &  &  &  & \multicolumn{1}{|c}{} &  &
  \multicolumn{1}{c:}{} &  &  & \zerobox{x} & \multicolumn{1}{c|}{\bu} &  &  &  &
   &  &  &  &  &  &  &  &
   &  & & \multicolumn{1}{|c}{}  &  &  & \multicolumn{1}{c|}{} &  &
   &  &  &  &  &  &  &  &
   &  &  &  &  &  &  &  \\
   &  &  &  &  &  &  &  &
   &  &  &  &  &  & \multicolumn{1}{|c}{} &  &
  \multicolumn{1}{c:}{} &  &  &  & \multicolumn{1}{c|}{} &  &  &  &
   &  &  &  &  &  &  &  &
   &  & & \multicolumn{1}{|c}{}  &  &  & \multicolumn{1}{c|}{} &  &
   &  &  &  &  &  &  &  &
   &  &  &  &  &  &  &  \\
   &  &  &  &  &  &  &  &
   &  &  &  &  &  & \multicolumn{1}{|c}{} &  &
  \multicolumn{1}{c:}{} &  &  &  & \multicolumn{1}{c|}{} &  &  &  &
   &  &  &  &  &  &  &  &
   &  & & \multicolumn{1}{|c}{} &  &  & \multicolumn{1}{c|}{} &  &
   &  &  &  &  &  &  &  &
   &  &  &  &  &  &  &  \\
   &  &  &  &  &  &  &  &
   &  &  &  &  &  & \multicolumn{1}{|c}{} &  &
  \multicolumn{1}{c:}{} &  &  & \zerobox[-3.6pt]{\bar{A}_{\mathrm{NW}}} & \multicolumn{1}{c|}{} &  &  &  &
   &  &  &  &  &  &  &  &
   &  & & \multicolumn{1}{|c}{} & \zerobox[3pt]{\bar{A}_{\mathrm{NE}}} &  & \multicolumn{1}{c|}{} &  &
   &  &  &  &  &  &  &  &
   &  &  &  &  &  &  &  \\
   &  &  &  &  &  &  &  &
   &  &  &  &  &  & \multicolumn{1}{|c}{} &  &
  \multicolumn{1}{c:}{} &  &  &  & \multicolumn{1}{c|}{} &  &  &  &
   &  &  &  &  &  &  &  &
   &  & & \multicolumn{1}{|c}{} &  &  & \multicolumn{1}{c|}{} &  &
   &  &  &  &  &  & \zerobox{z} &  &
   &  &  &  &  &  &  &  \\
\cline{15-21}\cline{36-39}
   &  &  &  &  &  &  &  &
   &  &  &  &  &  &  &  &
  \multicolumn{1}{c:}{} &  &  &  &  &  &  &  &
   &  &  &  &  &  &  &  &
   &  &  &  &  &  &  &  &
   &  &  &  &  & \bu &  &  &
   &  &  &  &  &  &  &  \\
\cline{15-21}\cline{36-39}
   &  &  &  &  &  &  &  &
   &  &  &  &  &  & \multicolumn{1}{|c}{} &  &
  \multicolumn{1}{c:}{} &  &  &  & \multicolumn{1}{c|}{} &  &  &  &
   &  &  &  &  &  &  &  &
   &  &  & \multicolumn{1}{|c}{}  &  &  & \multicolumn{1}{c|}{} &  &
   &  &  &  &  &  &  &  &
   &  &  &  &  &  &  &  \\
   &  &  &  &  &  &  &  &
   &  &  &  &  &  & \multicolumn{1}{|c}{} &  &
  \multicolumn{1}{c:}{} &  &  &  & \multicolumn{1}{c|}{} &  &  &  &
   &  &  &  &  &  &  &  &
   &  & & \multicolumn{1}{|c}{}  &  &  & \multicolumn{1}{c|}{} &  &
   &  &  &  &  &  &  &  &
   &  &  &  &  &  &  &  \\
\cline{1-12}
  \multicolumn{1}{|c}{} &  &  &  &  &  &  &  &
   &  &  & \multicolumn{1}{c|}{} &  &  & \multicolumn{1}{|c}{} &  &
  \multicolumn{1}{c:}{} &  &  & \zerobox[-3.6pt]{\bar{A}_{\mathrm{SW}}} & \multicolumn{1}{c|}{} &  &  &  &
   &  &  &  &  &  &  &  &
   &  & & \multicolumn{1}{|c}{}  & \zerobox[3pt]{\bar{A}_{\mathrm{SE}}} &  & \multicolumn{1}{c|}{} &  &
   &  &  &  &  &  &  &  &
   &  &  &  &  &  &  &  \\
  \multicolumn{1}{|c}{} &  &  &  &  &  &  &  &
   &  &  & \multicolumn{1}{c|}{} &  &  & \multicolumn{1}{|c}{} &  &
  \multicolumn{1}{c:}{} &  &  &  & \multicolumn{1}{c|}{} &  &  &  &
   &  &  &  &  &  &  &  &
   &  & & \multicolumn{1}{|c}{}  &  &  & \multicolumn{1}{c|}{} &  &
   &  &  &  &  &  &  &  &
   &  &  &  &  &  &  &  \\
\cline{36-45}\cline{47-51}
  \multicolumn{1}{|c}{} &  &  &  &  &  &  & \zerobox{A'_{\mathrm{S}}} &
   &  &  & \multicolumn{1}{c|}{} &  &  & \multicolumn{1}{|c}{} &  &
  \multicolumn{1}{c:}{} &  &  &  & \multicolumn{1}{c|}{} &  &  &  &
   &  &  &  &  &  &  &  &
   &  &  &  &  &  &  & \multicolumn{1}{|c}{} &
   &  &  &  & \multicolumn{1}{c|}{} &  & \multicolumn{1}{|c}{} &  &
   &  & \multicolumn{1}{c|}{} &  &  &  &  &  \\
\cline{15-21}
  \multicolumn{1}{|c}{} &  &  &  &  &  &  &  &
   &  &  & \multicolumn{1}{c|}{} &  &  &  &  &
  \multicolumn{1}{c:}{} &  &  &  &  &  &  &  &
   &  &  &  &  &  &  &  &
   &  &  &  &  &  &  & \multicolumn{1}{|c}{} &
   &  &  &  & \multicolumn{1}{c|}{} &  & \multicolumn{1}{|c}{} &  &
   &  & \multicolumn{1}{c|}{} &  &  &  &  &  \\
  \multicolumn{1}{|c}{} &  &  &  &  &  &  &  &
   &  &  & \multicolumn{1}{c|}{} &  &  &  &  &
  \multicolumn{1}{c:}{} &  &  &  &  &  &  &  &
   &  &  &  &  &  &  &  &
   &  &  &  &  &  &  & \multicolumn{1}{|c}{} &
   &  &  &  & \multicolumn{1}{c|}{} &  & \multicolumn{1}{|c}{} &  &
   &  & \multicolumn{1}{c|}{} &  &  &  &  &  \\
\cline{1-12}
   &  &  &  &  &  &  &  &
   &  &  &  &  &  &  &  &
  \multicolumn{1}{c:}{} &  &  &  &  &  &  &  &
   &  &  &  &  &  &  &  &
   &  &  &  &  &  &  & \multicolumn{1}{|c}{} &
   &  &  &  & \multicolumn{1}{c|}{} &  & \multicolumn{1}{|c}{} &  &
   &  & \multicolumn{1}{c|}{} &  &  &  &  &  \\
\cline{25-30}\cline{32-36}
   &  &  &  &  &  &  &  &
   &  &  &  &  &  &  &  &
  \multicolumn{1}{c:}{} &  &  &  &  &  &  &  &
  \multicolumn{1}{|c}{} &  &  &  &  & \multicolumn{1}{c|}{\bu} &  & \multicolumn{1}{|c}{} &
   &  &  & \multicolumn{1}{c|}{} &  &  &  & \multicolumn{1}{|c}{} &
   & \zerobox{B_{\mathrm{W}}} &  &  & \multicolumn{1}{c|}{} &  & \multicolumn{1}{|c}{} & \zerobox{B_{\mathrm{E}}} &
   &  & \multicolumn{1}{c|}{} &  &  &  &  &  \\
   &  &  &  &  &  &  &  &
   &  &  &  &  &  &  &  &
  \multicolumn{1}{c:}{} &  &  &  &  &  &  &  &
  \multicolumn{1}{|c}{} &  & \zerobox{\bar{B}_{\mathrm{SW}}} &  &  & \multicolumn{1}{c|}{\zerobox{y}} &  & \multicolumn{1}{|c}{} &
   & \zerobox{\bar{B}_{\mathrm{SE}}} &  & \multicolumn{1}{c|}{} &  &  &  & \multicolumn{1}{|c}{} &
   &  &  &  & \multicolumn{1}{c|}{} &  & \multicolumn{1}{|c}{} &  &
   &  & \multicolumn{1}{c|}{\bu} &  &  &  &  &  \\
   &  &  &  &  &  &  &  &
   &  &  &  &  &  &  &  &
  \multicolumn{1}{c:}{} &  &  &  &  &  &  &  &
  \multicolumn{1}{|c}{} &  &  &  &  & \multicolumn{1}{c|}{} &  & \multicolumn{1}{|c}{} &
   &  &  & \multicolumn{1}{c|}{} &  &  &  & \multicolumn{1}{|c}{} &
   &  &  &  & \multicolumn{1}{c|}{} &  & \multicolumn{1}{|c}{} &  &
   &  & \multicolumn{1}{c|}{\bu} &  &  &  &  &  \\
\cdashline{15-53}[2pt/1.5pt]
   &  &  &  &  &  &  &  &
   &  &  &  &  &  &  &  &
  \multicolumn{1}{c:}{} &  &  &  &  &  &  &  &
  \multicolumn{1}{|c}{} &  &  &  &  & \multicolumn{1}{c|}{} &  & \multicolumn{1}{|c}{} &
   &  &  & \multicolumn{1}{c|}{} &  &  &  & \multicolumn{1}{|c}{} &
   &  &  &  & \multicolumn{1}{c|}{\bu} &  & \multicolumn{1}{|c}{} &  &
   &  & \multicolumn{1}{c|}{} &  &  &  &  &  \\
   &  &  &  &  &  &  &  &
   &  &  &  &  &  &  &  &
  \multicolumn{1}{c:}{} &  &  &  &  &  &  &  &
  \multicolumn{1}{|c}{} &  &  &  &  & \multicolumn{1}{c|}{} &  & \multicolumn{1}{|c}{} &
   &  &  & \multicolumn{1}{c|}{} &  &  &  & \multicolumn{1}{|c}{} &
   &  &  &  & \multicolumn{1}{c|}{} &  & \multicolumn{1}{|c}{} &  &
   &  & \multicolumn{1}{c|}{} &  &  &  &  &  \\
\cline{25-30}\cline{32-36}
   &  &  &  &  &  &  &  &
   &  &  &  &  &  &  &  &
  \multicolumn{1}{c:}{} &  &  &  &  &  &  &  &
   &  &  &  &  &  &  &  &
   &  &  &  &  &  &  & \multicolumn{1}{|c}{} &
   &  &  &  & \multicolumn{1}{c|}{} &  & \multicolumn{1}{|c}{} &  &
   &  & \multicolumn{1}{c|}{} &  &  &  &  &  \\
   &  &  &  &  &  &  &  &
   &  &  &  &  &  &  &  &
   &  &  &  &  &  &  &  &
   &  &  &  &  &  &  &  &
   &  &  &  &  &  &  & \multicolumn{1}{|c}{} &
   &  &  &  & \multicolumn{1}{c|}{} &  & \multicolumn{1}{|c}{} &  &
   &  & \multicolumn{1}{c|}{} &  &  &  &  &  \\
   &  &  &  &  &  &  &  &
   &  &  &  &  &  &  &  &
   &  &  &  &  &  &  &  &
   &  &  &  &  &  &  &  &
   &  &  &  &  &  &  & \multicolumn{1}{|c}{} &
   &  &  &  & \multicolumn{1}{c|}{} &  & \multicolumn{1}{|c}{} &  &
   &  & \multicolumn{1}{c|}{} &  &  &  &  &  \\
\cline{23-28}\cline{40-45}\cline{47-56}
   &  &  &  &  &  &  &  &
   &  &  &  &  &  &  &  &
   &  &  &  &  &  & \multicolumn{1}{|c}{} &  &
   &  &  & \multicolumn{1}{c|}{} &  &  &  &  &
   &  &  &  &  &  &  &  &
   &  &  &  &  &  &  &  &
   &  &  & \multicolumn{1}{|c}{} &  &  &  & \multicolumn{1}{c|}{} \\
   &  &  &  &  &  &  &  &
   &  &  &  &  &  &  &  &
   &  &  &  &  &  & \multicolumn{1}{|c}{} &  &
   &  &  & \multicolumn{1}{c|}{} &  &  &  &  &
   &  &  &  &  &  &  &  &
   &  &  &  &  &  &  &  &
   &  &  & \multicolumn{1}{|c}{} &  &  &  & \multicolumn{1}{c|}{} \\
   &  &  &  &  &  &  &  &
   &  &  &  &  &  &  &  &
   &  &  &  &  &  & \multicolumn{1}{|c}{} &  &
   &  &  & \multicolumn{1}{c|}{} &  &  &  &  &
   &  &  &  &  &  &  &  &
   &  &  &  &  &  &  &  &
   &  &  & \multicolumn{1}{|c}{} &  &  &  & \multicolumn{1}{c|}{} \\
   &  &  &  &  &  &  &  &
   &  &  &  &  &  &  &  &
   &  &  &  &  &  & \multicolumn{1}{|c}{} &  &
   &  &  & \multicolumn{1}{c|}{} &  &  &  &  &
   &  &  &  &  &  &  &  &
   &  &  &  &  &  &  &  &
   &  &  & \multicolumn{1}{|c}{} &  &  &  & \multicolumn{1}{c|}{} \\
   &  &  &  &  &  &  &  &
   &  &  &  &  &  &  &  &
   &  &  &  &  &  & \multicolumn{1}{|c}{} &  &
   &  & \zerobox{B'_{\mathrm{W}}} & \multicolumn{1}{c|}{} &  &  &  &  &
   &  &  &  &  &  &  &  &
   &  &  &  &  &  &  &  &
   &  &  & \multicolumn{1}{|c}{} & \zerobox[1pt]{B'_{\mathrm{E}}} &  &  & \multicolumn{1}{c|}{} \\
   &  &  &  &  &  &  &  &
   &  &  &  &  &  &  &  &
   &  &  &  &  &  & \multicolumn{1}{|c}{} &  &
   &  &  & \multicolumn{1}{c|}{} &  &  &  &  &
   &  &  &  &  &  &  &  &
   &  &  &  &  &  &  &  &
   &  &  & \multicolumn{1}{|c}{} &  &  &  & \multicolumn{1}{c|}{\bu} \\
   &  &  &  &  &  &  &  &
   &  &  &  &  &  &  &  &
   &  &  &  &  &  & \multicolumn{1}{|c}{} &  &
   &  &  & \multicolumn{1}{c|}{} &  &  &  &  &
   &  &  &  &  &  &  &  &
   &  &  &  &  &  &  &  &
   &  &  & \multicolumn{1}{|c}{} &  &  &  & \multicolumn{1}{c|}{\bu} \\
   &  &  &  &  &  &  &  &
   &  &  &  &  &  &  &  &
   &  &  &  &  &  & \multicolumn{1}{|c}{} &  &
   &  &  & \multicolumn{1}{c|}{\bu} &  &  &  &  &
   &  &  &  &  &  &  &  &
   &  &  &  &  &  &  &  &
   &  &  & \multicolumn{1}{|c}{} &  &  &  & \multicolumn{1}{c|}{} \\
   &  &  &  &  &  &  &  &
   &  &  &  &  &  &  &  &
   &  &  &  &  &  & \multicolumn{1}{|c}{} &  &
   &  &  & \multicolumn{1}{c|}{} &  &  &  &  &
   &  &  &  &  &  &  &  &
   &  &  &  &  &  &  &  &
   &  &  & \multicolumn{1}{|c}{} &  &  &  & \multicolumn{1}{c|}{} \\
   &  &  &  &  &  &  &  &
   &  &  &  &  &  &  &  &
   &  &  &  &  &  & \multicolumn{1}{|c}{} &  &
   &  &  & \multicolumn{1}{c|}{} &  &  &  &  &
   &  &  &  &  &  &  &  &
   &  &  &  &  &  &  &  &
   &  &  & \multicolumn{1}{|c}{} &  &  &  & \multicolumn{1}{c|}{} \\
   &  &  &  &  &  &  &  &
   &  &  &  &  &  &  &  &
   &  &  &  &  &  & \multicolumn{1}{|c}{} &  &
   &  &  & \multicolumn{1}{c|}{} &  &  &  &  &
   &  &  &  &  &  &  &  &
   &  &  &  &  &  &  &  &
   &  &  & \multicolumn{1}{|c}{} &  &  &  & \multicolumn{1}{c|}{} \\
   &  &  &  &  &  &  &  &
   &  &  &  &  &  &  &  &
   &  &  &  &  &  & \multicolumn{1}{|c}{} &  &
   &  &  & \multicolumn{1}{c|}{} &  &  &  &  &
   &  &  &  &  &  &  &  &
   &  &  &  &  &  &  &  &
   &  &  & \multicolumn{1}{|c}{} &  &  &  & \multicolumn{1}{c|}{} \\
\cline{23-28}\cline{52-56}
\end{array}
\)%

%% file: references.bib
@string{JACM = {J.~ACM}}

@string{SIDMA = {SIAM J.~Discrete Mathematics}}

@string{FOCS = {Proc.~FOCS}}

@string{JCTA = {J.~Combin. Theory Ser.~A}}

@string{SODA = {Proc.~SODA}}

@string{DM = {Discrete Mathematics}}

@article{Cibulka09,
  author       = {Josef Cibulka},
  title        = {On constants in the {F}{\"{u}}redi-{H}ajnal and the {S}tanley-{W}ilf
                  conjecture},
  journal      = {J. Comb. Theory, Ser. {A}},
  volume       = {116},
  number       = {2},
  pages        = {290--302},
  year         = {2009},
  url          = {https://doi.org/10.1016/j.jcta.2008.06.003},
  doi          = {10.1016/j.jcta.2008.06.003},
}

@inproceedings{DemaineHIKP09,
  author =    "Erik D. Demaine and Dion Harmon and John Iacono and Daniel M. Kane and Mihai P{\v a}tra{\c s}cu",
  title =     "The Geometry of Binary Search Trees",
  booktitle = "Proc. 20th ACM-SIAM Symposium on Discrete Algorithms (SODA)",
  year =      "2009",
  pages =		"496--505",
}

@article{Fulek09,
  author =	{Radoslav Fulek},
  title = {Linear bound on extremal functions of some forbidden patterns in 0--1 matrices},
  journal = DM,
  volume =	{309},
  issue =		{6},
  year =		{2009},
  pages =		{1736--1739},
}

@incollection{FurediS13,
    author= {F{\"u}redi, Zolt{\'a}n and Simonovits, Mikl{\'o}s},
    editor= {Lov{\'a}sz, L{\'a}szl{\'o} and Ruzsa, Imre Z.
                and S{\'o}s, Vera T.},
    title={The History of Degenerate (Bipartite) Extremal Graph Problems},
    booktitle={Erd{\H{o}}s Centennial},
    series ={Bolyai Society Mathematical Studies},
    number ={25},
    year={2013},
    publisher={Springer},
    Xaddress={Berlin, Heidelberg},
    pages={169--264},
    doi={10.1007/978-3-642-39286-3_7},
    url={https://doi.org/10.1007/978-3-642-39286-3_7},
}

@article{Furedi90,
  author    = {Zolt\'{a}n F{\"{u}}redi},
  title     = {The maximum number of unit distances in a convex $N$-gon},
  journal   = JCTA,
  volume    = {55},
  number    = {2},
  year      = {1990},
  pages     = {316--320},
}

@article{FurediH92,
  author    = {Zolt\'{a}n F{\"{u}}redi and P\'{e}ter Hajnal},
  title     = {Davenport-{S}chinzel theory of matrices},
  journal   = DM,
  volume    = {103},
  number    = {3},
  year      = {1992},
  pages     = {233--251},
}

@article{Geneson09,
  title	 = {Extremal functions of forbidden double permutation matrices},
  journal	= JCTA,
  volume	= {116},
  year 	= {2009},
  number	= {7},
  pages 	= {1235--1244},
  author	= {Jesse T. Geneson},
}

@article {HS86,
    AUTHOR = {Hart, Sergiu and Sharir, Micha},
     TITLE = {Nonlinearity of {D}avenport-{S}chinzel sequences and of
              generalized path compression schemes},
   JOURNAL = {Combinatorica},
  FJOURNAL = {Combinatorica. An International Journal of the J\'anos Bolyai
              Mathematical Society},
    VOLUME = {6},
      YEAR = {1986},
    NUMBER = {2},
     PAGES = {151--177},
}

@article{Keszegh09,
  author = {Bal\'{a}zs Keszegh},
  title = {On linear forbidden submatrices},
  journal = JCTA,
  year = {2009},
  volume = {116},
  number = {1},
  pages = {232--241},
}

@incollection {Klazar00,
    AUTHOR = {Martin Klazar},
     TITLE = {The {F}\"uredi-{H}ajnal conjecture implies the {S}tanley-{W}ilf conjecture},
 BOOKTITLE = {Formal power series and algebraic combinatorics ({M}oscow, 2000)},
     PAGES = {250--255},
 PUBLISHER = {Springer},
   ADDRESS = {Berlin},
      YEAR = {2000},
}

@article {Klazar92,
    AUTHOR = {Martin Klazar},
     TITLE = {A general upper bound in extremal theory of sequences},
   JOURNAL = {Comment. Math. Univ. Carolin.},
  FJOURNAL = {Commentationes Mathematicae Universitatis Carolinae},
    VOLUME = {33},
      YEAR = {1992},
    NUMBER = {4},
     PAGES = {737--746},
}

@article{KovariST54,
    AUTHOR = {K\H{o}vari, Tam\'{a}s and S{\'o}s, Vera~T. and Tur{\'a}n, P\'{a}l},
     TITLE = {On a problem of {K}. {Z}arankiewicz},
   JOURNAL = {Colloquium Math.},
    VOLUME = {3},
      YEAR = {1954},
     PAGES = {50--57},
}

@article {MarcusT04,
    AUTHOR = {Marcus, Adam and Tardos, G\'{a}bor},
     TITLE = {Excluded permutation matrices and the {S}tanley-{W}ilf conjecture},
   JOURNAL = JCTA,
  FJOURNAL = {Journal of Combinatorial Theory. Series A},
    VOLUME = {107},
      YEAR = {2004},
    NUMBER = {1},
     PAGES = {153--160},
}

@article{Nivasch10,
	author = {Gabriel Nivasch},
	title = {Improved bounds and new techniques for {D}avenport-{S}chinzel sequences and their generalizations},
	journal = JACM,
	volume = {57},
	number = {3},
	articleno = {17},
	year = {2010},
}

@article {PachTardos06,
    AUTHOR = {Pach, J\'{a}nos and Tardos, G\'{a}bor},
     TITLE = {Forbidden paths and cycles in ordered graphs and matrices},
   JOURNAL = {Israel J. Math.},
  FJOURNAL = {Israel Journal of Mathematics},
    VOLUME = {155},
      YEAR = {2006},
     PAGES = {359--380},
}

@inproceedings{Pettie10a,
  author      = {Seth Pettie},
  title       = {Applications of forbidden 0--1 matrices to search tree- and path compression-based data structures},
  booktitle   = {Proceedings 21st ACM-SIAM Symposium on Discrete Algorithms (SODA)},
  year        = {2010},
  pages       = {1457--1467}
}

@article{Pettie-FH11,
  author	= {Seth Pettie},
  title		= {Degrees of nonlinearity in forbidden 0--1 matrix problems},
  journal	= {Discrete Mathematics},
  year	= {2011},
  volume	= {311},
  pages	= {2396--2410},
}

@article{Pettie-GenDS11,
  author	= {Seth Pettie},
  title		= {Generalized {D}avenport-{S}chinzel sequences and their 0--1 matrix counterparts},
  journal	= {J. Comb. Theory Ser. A},
  year	= {2011},
  volume	= {118},
  number	= {6},
  pages	= {1863--1895},
}

@inproceedings{Pettie-SoCG11,
  author	= {Seth Pettie},
  title		= {On the structure and composition of forbidden sequences, with geometric applications},
  booktitle	= {Proceedings 27th Annual Symposium on Computational Geometry (SoCG)},
  year	= {2011},
  pages	= {370--379},
}

@article{Pettie15-SIDMA,
  author    = {Seth Pettie},
  title     = {Three Generalizations of {D}avenport-{S}chinzel Sequences},
  journal   = SIDMA,
  volume    = {29},
  number    = {4},
  pages     = {2189--2238},
  year      = {2015},
  doi       = {10.1137/140968574},
}

@article{Reiman58,
	author = 	{I. Reiman},
	title = 	{\"{U}ber ein {P}roblem von {K}. {Z}arankiewicz},
	journal = 	{Acta. Math. Acad. Sci. Hungary},
	volume = 	{9},
	pages = 	{269--273},
	year = 	{1958},
}

@article {ST85,
    AUTHOR = {Sleator, Daniel~D. and Tarjan, Robert~E.},
     TITLE = {Self-adjusting binary search trees},
   JOURNAL = JACM,
  FJOURNAL = {Journal of the Association for Computing Machinery},
    VOLUME = {32},
      YEAR = {1985},
    NUMBER = {3},
     PAGES = {652--686},
}

@article {Tardos05,
    AUTHOR = {Tardos, G\'{a}bor},
     TITLE = {On 0--1 matrices and small excluded submatrices},
   JOURNAL = JCTA,
  FJOURNAL = {Journal of Combinatorial Theory. Series A},
    VOLUME = {111},
      YEAR = {2005},
    NUMBER = {2},
     PAGES = {266--288},
}

@misc{Fox13,
      title={{S}tanley-{W}ilf limits are typically exponential}, 
      author={Jacob Fox},
      year={2013},
      note={arXiv:1310.8378},
      eprint={1310.8378},
      archivePrefix={arXiv},
      primaryClass={math.CO},
      url={https://arxiv.org/abs/1310.8378}, 
}

@inproceedings{CibulkaK17,
  author       = {Josef Cibulka and
                  Jan Kyn\v{c}l},
  title        = {Better upper bounds on the {F}{\"{u}}redi-{H}ajnal limits of permutations},
  booktitle    = {Proceedings of the 28th {ACM-SIAM} Symposium on Discrete
                  Algorithms (SODA)},
  pages        = {2280--2293},
  year         = {2017},
  url          = {https://doi.org/10.1137/1.9781611974782.150},
  doi          = {10.1137/1.9781611974782.150},
}

@inproceedings{ChalermsookGKMS15,
  author       = {Parinya Chalermsook and
                  Mayank Goswami and
                  L{\'{a}}szl{\'{o}} Kozma and
                  Kurt Mehlhorn and
                  Thatchaphol Saranurak},
  title        = {Pattern-Avoiding Access in Binary Search Trees},
  booktitle    = {Proceedings of the 56th {IEEE} 
                    Symposium on Foundations of Computer Science ({FOCS})},
  pages        = {410--423},
  year         = {2015},
  url          = {https://doi.org/10.1109/FOCS.2015.32},
  doi          = {10.1109/FOCS.2015.32},
}

@inproceedings{ChalermsookGJAPY23,
  author       = {Parinya Chalermsook and
                  Manoj Gupta and
                  Wanchote Jiamjitrak and
                  Nidia Obscura Acosta and
                  Akash Pareek and
                  Sorrachai Yingchareonthawornchai},
  title        = {Improved Pattern-Avoidance Bounds for Greedy {BST}s via Matrix Decomposition},
  booktitle    = {Proceedings of the 34th {ACM-SIAM} Symposium on Discrete Algorithms (SODA)},
  pages        = {509--534},
  year         = {2023},
  url          = {https://doi.org/10.1137/1.9781611977554.ch22},
  doi          = {10.1137/1.9781611977554.ch22},
}

@inproceedings{ChalermsookGKMS15b,
  author       = {Parinya Chalermsook and
                  Mayank Goswami and
                  L{\'{a}}szl{\'{o}} Kozma and
                  Kurt Mehlhorn and
                  Thatchaphol Saranurak},
  title        = {Greedy Is an Almost Optimal Deque},
  booktitle    = {Proceedings 14th Symposium on Algorithms and Data Structures},
  series       = {Lecture Notes in Computer Science},
  volume       = {9214},
  pages        = {152--165},
  Xpublisher    = {Springer},
  year         = {2015},
  url          = {https://doi.org/10.1007/978-3-319-21840-3\_13},
  doi          = {10.1007/978-3-319-21840-3\_13},
}

@article{KorandiTTW19,
  author       = {D{\'{a}}niel Kor{\'{a}}ndi and
                  G{\'{a}}bor Tardos and
                  Istv{\'{a}}n Tomon and
                  Craig Weidert},
  title        = {On the {T}ur{\'{a}}n number of ordered forests},
  journal      = {J. Comb. Theory, Ser. {A}},
  volume       = {165},
  pages        = {32--43},
  year         = {2019},
  url          = {https://doi.org/10.1016/j.jcta.2019.01.006},
  doi          = {10.1016/j.jcta.2019.01.006},
}

@inproceedings{GuillemotM14,
  author       = {Sylvain Guillemot and
                  D{\'{a}}niel Marx},
  title        = {Finding small patterns in permutations in linear time},
  booktitle    = {Proceedings 25th Annual {ACM-SIAM} Symposium on Discrete Algorithms ({SODA})},
  pages        = {82--101},
  year         = {2014},
  url          = {https://doi.org/10.1137/1.9781611973402.7},
  doi          = {10.1137/1.9781611973402.7},
}

@inproceedings{ChalermsookPY24,
  author       = {Parinya Chalermsook and
                  Seth Pettie and
                  Sorrachai Yingchareonthawornchai},
  title        = {Sorting Pattern-Avoiding Permutations via 0-1 Matrices Forbidding
                  Product Patterns},
  booktitle    = {Proceedings 35th {ACM-SIAM} Symposium on Discrete Algorithms (SODA)},
  pages        = {133--149},
  year         = {2024},
  url          = {https://doi.org/10.1137/1.9781611977912.7},
  doi          = {10.1137/1.9781611977912.7},
}

@inproceedings{PettieT24,
  author       = {Seth Pettie and
                  G{\'{a}}bor Tardos},
  title        = {On the Extremal Functions of Acyclic Forbidden 0-1 Matrices},
  booktitle    = {Proceedings 35th {ACM-SIAM} Symposium on Discrete Algorithms ({SODA})},
  pages        = {1166--1176},
  year         = {2024},
  url          = {https://doi.org/10.1137/1.9781611977912.45},
  doi          = {10.1137/1.9781611977912.45},
}

@article{PettieT25,
  author       = {Seth Pettie and
                  G{\'{a}}bor Tardos},
  title        = {A Refutation of the {P}ach--{T}ardos Conjecture for 0--1 Matrices},
  journal      = {Comb.},
  volume       = {45},
  number       = {6},
  pages        = {59},
  year         = {2025},
  url          = {https://doi.org/10.1007/s00493-025-00187-7},
  doi          = {10.1007/S00493-025-00187-7},
}

@article{Geneson26,
  author       = {Jesse Geneson},
  title        = {Linear extremal bounds for a family of forbidden 0--1 matrices},
  journal      = {CoRR},
  volume       = {abs/2607.16463},
  year         = {2026},
  url          = {https://doi.org/10.48550/arXiv.2607.16463},
  doi          = {10.48550/ARXIV.2607.16463},
  eprinttype   = {arXiv},
  eprint       = {2607.16463},
}

@article{Geneson26b,
  author       = {Jesse Geneson},
  title        = {Typical growth of the {F}\"uredi-{H}ajnal and {S}tanley-{W}ilf limits},
  journal      = {CoRR},
  volume       = {abs/2610.02707},
  year         = {2026},
  url          = {https://doi.org/10.48550/arXiv.2610.02707},
  doi          = {10.48550/ARXIV.2610.02707},
  eprinttype   = {arXiv},
  eprint       = {2610.02707},
}

@article{GenesonTsai20,
  author       = {Jesse Geneson and
                  Shen{-}Fu Tsai},
  title        = {Sharper Bounds and Structural Results for Minimally Nonlinear 0--1
                  Matrices},
  journal      = {Electron. J. Comb.},
  volume       = {27},
  number       = {4},
  pages        = {P4.24},
  year         = {2020},
  url          = {https://doi.org/10.37236/7801},
  doi          = {10.37236/7801},
}

@article{GenesonShen,
  author       = {Jesse Geneson and
                  Lilly Shen},
  title        = {Linear bounds on matrix extremal functions using visibility hypergraphs},
  journal      = {Discret. Math.},
  volume       = {338},
  number       = {12},
  pages        = {2437--2441},
  year         = {2015},
  url          = {https://doi.org/10.1016/j.disc.2015.06.017},
  doi          = {10.1016/J.DISC.2015.06.017},
}

@article{GenesonT17,
  author       = {Jesse T. Geneson and
                  Peter M. Tian},
  title        = {Extremal functions of forbidden multidimensional matrices},
  journal      = {Discret. Math.},
  volume       = {340},
  number       = {12},
  pages        = {2769--2781},
  year         = {2017},
  url          = {https://doi.org/10.1016/j.disc.2017.08.017},
  doi          = {10.1016/J.DISC.2017.08.017},
}

@article{CrowdMath2018, 
    title={Bounds on Parameters of Minimally Nonlinear Patterns}, 
    volume={25}, 
    url={https://www.combinatorics.org/ojs/index.php/eljc/article/view/v25i1p5}, 
    DOI={10.37236/6735}, 
    number={1}, 
    journal={The Electronic Journal of Combinatorics}, author={CrowdMath, P. A.}, 
    year={2018},
    pages={P1.5}, 
}

@article{BonnetNMST24,
  author       = {{\'{E}}douard Bonnet and
                  Jaroslav Nesetril and
                  Patrice Ossona de Mendez and
                  Sebastian Siebertz and
                  St{\'{e}}phan Thomass{\'{e}}},
  title        = {Twin-width and permutations},
  journal      = {Log. Methods Comput. Sci.},
  volume       = {20},
  number       = {3},
  year         = {2024},
  url          = {https://doi.org/10.46298/lmcs-20(3:4)2024},
  doi          = {10.46298/LMCS-20(3:4)2024},
}

@article{BonnetGMSTT24,
  author       = {{\'{E}}douard Bonnet and
                  Ugo Giocanti and
                  Patrice Ossona de Mendez and
                  Pierre Simon and
                  St{\'{e}}phan Thomass{\'{e}} and
                  Szymon Torunczyk},
  title        = {Twin-Width {IV:} Ordered Graphs and Matrices},
  journal      = {J. {ACM}},
  volume       = {71},
  number       = {3},
  pages        = {21},
  year         = {2024},
  url          = {https://doi.org/10.1145/3651151},
  doi          = {10.1145/3651151},
}

@article{GishbolinerL26,
  author    = {Lior Gishboliner and Xiangyu Li},
  title     = {Extremal numbers of acyclic matrices},
  journal   = {CoRR},
  volume    = {abs/2609.20726},
  year      = {2026},
  url       = {http://arxiv.org/abs/2609.20726},
  archivePrefix = {arXiv},
  eprint    = {2609.20726},
}

@article{AppelHaken1977,
  author    = {Appel, Kenneth and Haken, Wolfgang},
  title     = {Every planar map is four colorable. {Part {I}: Discharging}},
  journal   = {Illinois Journal of Mathematics},
  volume    = {21},
  number    = {3},
  pages     = {429--490},
  year      = {1977},
  url       = {http://projecteuclid.org}
}

@article{AppelHakenKoch1977,
  author    = {Appel, Kenneth and Haken, Wolfgang and Koch, John},
  title     = {Every planar map is four colorable. {Part {II}: Reducibility}},
  journal   = {Illinois Journal of Mathematics},
  volume    = {21},
  number    = {3},
  pages     = {491--567},
  year      = {1977},
}

@inproceedings{InoueKMMTT26,
  author    = {Yuta Inoue and {Ken-ichi} Kawarabayashi and Atsuyuki Miyashita and Bojan Mohar and Carsten Thomassen and Mikkel Thorup},
  title     = {The Four Color Theorem with Linearly Many Reducible Configurations and Near-Linear Time Coloring},
  booktitle = {Proceedings of the 67th IEEE Annual Symposium on Foundations of Computer Science (FOCS)},
  year      = {2026},
  note      = {Best Paper Award}
}

@misc{InoueKMMTT26b,
      title={Four-Coloring Planar Graphs in Linear Time}, 
      author={Yuta Inoue and {Ken-ichi} Kawarabayashi and Atsuyuki Miyashita and Bojan Mohar and Carsten Thomassen and Mikkel Thorup},
      year={2026},
      eprint={2610.39753},
      archivePrefix={arXiv},
      primaryClass={cs.DS},
      url={https://arxiv.org/abs/2610.39753}, 
}

@article{DuchetHVM83,
  author       = {Pierre Duchet and
                  Yahya Ould Hamidoune and
                  Michel Las Vergnas and
                  Henry Meyniel},
  title        = {Representing a planar graph by vertical lines joining different levels},
  journal      = {Discret. Math.},
  volume       = {46},
  number       = {3},
  pages        = {319--321},
  year         = {1983},
  url          = {https://doi.org/10.1016/0012-365X(83)90128-0},
  doi          = {10.1016/0012-365X(83)90128-0},
}

@article{TamassiaT86,
  author       = {Roberto Tamassia and
                  Ioannis G. Tollis},
  title        = {A Unified Approach to Visibility Representations of Planar Graphs},
  journal      = {Discret. Comput. Geom.},
  volume       = {1},
  pages        = {321--341},
  year         = {1986},
  url          = {https://doi.org/10.1007/BF02187705},
  doi          = {10.1007/BF02187705},
}
